\documentclass[12pt,reqno]{amsart} 
\usepackage{amsmath,amscd,amssymb,amsfonts,amsthm}
\usepackage[mathscr]{eucal}
\usepackage{hyperref}

\title[Berezinians, Exterior Powers and Recurrent Sequences]{Berezinians,
Exterior Powers and Recurrent Sequences}

\author{H.~M.~Khudaverdian}
\author{Th.~Th. Voronov}

\address{H. K., Th. V.\,: Department of Mathematics, University of
Manchester Institute of Science and Technology (UMIST), Manchester M60 1QD,
United Kingdom}

\email{theodore.voronov@umist.ac.uk, khudian@umist.ac.uk}

\address{H. K.\,: G.~S.~Sahakian~Department~of~Theoretical ~Physics,
Yerevan State University, 1 A. Manoukian Street, 375049 Yerevan, Armenia}

\address{H. K.\,: Joint Institute for Nuclear Research, Dubna 141980,
Russia}

\numberwithin{equation}{section} 

\def\co{\colon\thinspace}
\newcommand{\RHS}{{r.h.s.}}

\newcommand{\mathcall}{\EuScript}

\renewcommand{\leq}{\leqslant}
\renewcommand{\geq}{\geqslant}

\DeclareMathOperator{\Ber}{Ber} \DeclareMathOperator{\Berp}{Ber^+}
\DeclareMathOperator{\Berm}{Ber^{\displaystyle --}}
\DeclareMathOperator{\Res}{Res} 
\DeclareMathOperator{\Tr}{Tr} \DeclareMathOperator{\GL}{GL}
\DeclareMathOperator{\adj}{adj} \DeclareMathOperator{\dsum}{{\oplus}}
\newcommand{\Rp}{R^+}
\newcommand{\Rm}{R^-}

\def \c  {\boldsymbol{c}}
\def\p{\partial}

\def \Z{{\mathbb Z}}
\def \ZZ{{\mathbb{Z}_2}}

\def \R {{R}}

\def \B   {{B}}
\def \P   {{\H}}

\renewcommand{\a}{{\alpha}}
\renewcommand{\b}{{\beta}}
\newcommand{\e}{{\varepsilon}}
\newcommand{\f}{{\varphi}}
\renewcommand{\l}{{\lambda}}
\newcommand{\m}{{\mu}}
\newcommand{\s}{{\sigma}}

\newcommand{\vt}{{\tilde v}}
\newcommand{\ut}{{\tilde u}}

\newcommand{\dif}{{\gamma}}
\newcommand{\G}{{\Gamma}}
\newcommand{\seq}{{\mathcall{S}}}
\newcommand{\D}{{\Delta}}

\newcommand{\be}{\boldsymbol{e}}
\newcommand{\bx}{\boldsymbol{x}}
\newcommand{\by}{\boldsymbol{y}}
\newcommand{\bv}{\boldsymbol{v}}
\renewcommand{\r}{\boldsymbol{\rho}}

\newcommand{\wed}{{\Lambda}}
\newcommand{\wz}{\wed_z}
\newcommand{\wzp}{\wed_{-z\Pi}}
\newcommand{\HC}{\mathcall{P}}
\renewcommand{\H}{\mathcall{H}}

\newtheorem{thm}{Theorem}
\newtheorem{prop}{Proposition}

\newtheorem*{coro}{Corollary}

\theoremstyle{definition}
\newtheorem*{de}{Definition}
\newtheorem*{exe}{Example}
\newtheorem{ex}{Example}[section]
\theoremstyle{remark}

\newtheorem*{Rem}{Remark}

\begin{document}

\begin{abstract}
We study power expansions of the characteristic function of  a linear operator
$A$ in a $p|q$-dimensional superspace $V$. We show that   traces of  exterior
powers of $A$ satisfy universal recurrence relations of period $q$.
`Underlying' recurrence relations hold in the Grothendieck ring of
representations of $\GL(V)$. They are expressed by vanishing of certain Hankel
determinants of order $q+1$ in this ring, which generalizes the vanishing of
sufficiently high exterior powers of an ordinary vector space. In particular,
this allows to explicitly express the Berezinian of an operator as a rational
function of traces. We  analyze the Cayley--Hamilton identity in a superspace.
Using the geometric meaning of the Berezinian we also give a simple formulation
of the analog of Cramer's rule.
\end{abstract}

\maketitle \tableofcontents

\section{Introduction}

\subsection{}
In this paper we study the Berezinians of linear operators in a superspace and
in particular the characteristic function $R_A(z)=\Ber(1+zA)$,
%
%
where $z$ is a complex variable. Our principal tool is the two power expansions
of $R_A(z)$, at zero and at infinity. We also study a similar rational function
taking values in a Grothendieck ring. The main results are as follows.

For an arbitrary even linear operator $A$ in a $p|q$-dimensional superspace $V$
we establish universal recurrence relations satisfied by the traces
$\Tr\wed^kA$ and $\Tr\Sigma^kA$ of the induced action in the exterior powers
$\wed^k (V)$ and the `dual exterior powers' $\Sigma^k(V)=\Ber V\otimes
\wed^{p-k}V^*$ (Theorem~\ref{thmrecurrence}, formulae~\eqref{eqhankelg} and
\eqref{eqhankelc}). We obtain similar fundamental recurrence relations
satisfied by the spaces $\wed^k (V)$ and $\Sigma^n(V)$ in a suitable
Grothendieck ring, and underlying the relations for traces
(Theorems~\ref{thmhankelgroth} and~\ref{thmgroth2}). In particular, we show how
$\Tr\Sigma^kA$, which are rational functions of $A$, can be obtained from the
polynomial invariants $\Tr\wed^kA$ by a sort of ``analytic continuation''. Our
considerations lead to effective formulae. For the Berezinian $\Ber A$ we
obtain an invariant explicit formula  expressing it as the ratio of two Hankel
determinants built of $\Tr\wed^k A$:
\begin{equation*}
    \Ber A=\frac{\begin{vmatrix}
      \Tr\wed^{p-q}A & \ldots & \Tr\wed^pA \\
      \ldots & \ldots & \ldots \\
      \Tr\wed^pA & \ldots & \Tr\wed^{p+q}A \\
    \end{vmatrix}}{\begin{vmatrix}
      \Tr\wed^{p-q+2}A & \ldots & \Tr\wed^{p+1}A \\
      \ldots & \ldots & \ldots \\
      \Tr\wed^{p+1}A & \ldots & \Tr\wed^{p+q}A \\
    \end{vmatrix}}\,.
\end{equation*}
One can relate these determinants with  characters of polynomial
representations of the general linear supergroup corresponding to particular
Young diagrams.

Besides this, we discuss two other related topics. For an analog of the
Cayley--Hamilton theorem, we analyze the problem of a minimal annihilating
polynomial of a linear operator in a superspace  and show how it can be
obtained from the characteristic function $R_A(z)$. It should be emphasized
that in the supercase the rational characteristic function $R_A(z)=\Ber(1+zA)$
is a more fundamental object than such a `characteristic polynomial', which can
be built from it. We also study an analog of Cramer's rule for the supercase
and give for it a geometric proof.

\subsection{Motivation and background.}  Recall that the Berezinian is the analog of the
determinant for the $\ZZ$-graded (= super) situation. It was discovered by
F.~A.~Berezin in  his studies of second quantization and   integration over odd
variables. See~\cite{berezin:secondre, berezin:antire} and references therein.
The main feature of $\Ber A$ is that it is not a polynomial in the matrix
entries, but a fraction. In the standard definition
\begin{equation*}
   \Ber A:=
   {\det \left(A_{00}-A_{01}A^{-1}_{11}A_{10}\right)}{(\det A_{11})^{-1}}\,,
\end{equation*}
where $A_{00}$, $A_{01}$, $A_{10}$ and $A_{11}$ are the matrix blocks of $A$,
the numerator and denominator do not have independent invariant meaning.
Exactly because $\Ber A$ is non-polynomial, integration theory in the supercase
is non-trivial. In particular, it is well known that the straightforward
generalization of the exterior algebra by standard tensor tools transferred to
the $\ZZ$-graded situation, is not sufficient, because it is not related with
the Berezinian and hence with integration over supermanifolds (for a survey
see, e.g., \cite{tv:git, tv:susy}). The simplest objects that one has to
consider besides the naive exterior powers $\wed^k(V)$ are the `dual exterior
powers' $\Sigma^k(V):=\Ber V\otimes \wed^{p-k}V^*$ introduced by Bernstein and
Leites~\cite{berl:int} (when $V$ is the space of covectors on a supermanifold
the elements of $\Sigma^k(V)$ are called integral forms).

As we show in this paper,  there are surprising ``hidden relations'' between
the naive exterior powers $\wed^k(V)$ and the Berezinian, so they are closer
than might be expected. This is seen by the comparing of the two expansions of
the characteristic function of a linear operator: the expansion at zero gives
the traces in $\wed^k(V)$, while the expansion at infinity gives the traces in
$\Sigma^k(V)$, including the Berezinian. Hence the relations between
$\wed^k(V)$ and $\Ber$ can be perceived as an `analytic continuation of a
rational function from a neighborhood of zero to the neighborhood of infinity'.
(There is an analogy with rational numbers: the ordinary decimal expansion
corresponds to an expansion near infinity, while a $p$-adic expansion
corresponds to an expansion at zero.) Formal analogs of these expansions yield
underlying relations in the Grothendieck ring.

Let us explain the position of these results in comparison with the familiar
picture of operators acting in purely even vector spaces. For a vector space
$V$ of dimension $n$ all exterior powers starting from $\wed^{n+1} (V)$,
vanish. Therefore all the traces $\Tr\wed^{k}A$, $k>n$, identically vanish.
Also, the top exterior power $\wed^n(V)$ is the same as the one-dimensional
space $\det V$, and this gives rise to natural isomorphisms  $\det V\otimes
\wed^{n-k}(V^*) \cong \wed^k(V)$ (`duality'). In the $\ZZ$-graded case, for a
vector space $V$ of dimension $p|q$, there is an infinite sequence of the
exterior powers $\wed^{k} (V)$, which does not terminate. Likewise,  there is
an infinite sequence of the spaces $\Sigma^k(V)=\Ber V\otimes \wed^{p-k}V^*$,
stretching to the left, which are now essentially different from $\wed^{k}
(V)$. In this paper we establish the following relations in the Grothendieck
ring:
\begin{equation*}
    \begin{vmatrix}
      \G_k  & \dots & \G_{k+q}  \\
      \dots & \dots & \dots \\
      \G_{k+q}  & \dots & \G_{k+2q}  \\
    \end{vmatrix}=0\,,
\end{equation*}
where $\G_k=\wed^k V-(-\Pi)^q\Sigma^{k+q} V$, for all $k\in \Z$. ($\Pi$ is the
parity shift functor.) Taken in the range of $k$ where both $\wed^kV$ and
$\Sigma^k V$ are not zero, it gives the proper replacement for the classical
`duality isomorphisms'. At the same time, its   corollary
\begin{equation*}
    \begin{vmatrix}
      \wed^kV & \dots & \wed^{k+q}V \\
      \dots & \dots & \dots \\
      \wed^{k+q}V & \dots & \wed^{k+2q}V \\
    \end{vmatrix}=0\,,
\end{equation*}
for  $k\geq p-q+1$ replaces the vanishing of the sufficiently high exterior
powers in the classical case.

The Cayley--Hamilton theorem is closely related with identities for traces. In
the classical case, $\chi_A(A)=0$ for the characteristic polynomial
$\chi_A(z)=\det(A-z)$ of a linear operator in an $n$-dimensional space, gives
relations for the powers of $A$. It can be deduced from the identity
$\Tr\wed^{n+1}A=0$ by varying it w.r.t. $A$, and, conversely, it implies
identities for traces. Now, in the $\ZZ$-graded case, of course, any even
operator  satisfies the same polynomial relation as in the classics with
$n=p+q$. The trouble, however, is how to give a meaning to the coefficients of
this relation as invariants of the operator. This has been a source of
confusion of many attempts to generalize the Cayley--Hamilton theorem to the
supercase that can be found in the literature. In this paper we explain how the
`naive' Cayley--Hamilton identity (if one forgets about the $\ZZ$-grading) and
an identity obtained by varying the relation for traces following from the
second formula above, give the same thing. The subordinate role of the
`Cayley--Hamilton polynomial' in the supercase as compared to the `true'
characteristic function $\Ber(A-z)$ (or the equivalent $R_A(z)=\Ber(1+zA)$), is
clearly seen.

Notice that in the last fifteen years there has been an active work on
non-commu\-ta\-tive generalizations of determinants initiated by Gelfand and
Retakh (see~\cite{gelfand:retakh91, gelfand:retakh02}), non-commutative Vieta
formulae~\cite{schwarz:vieta1, schwarz:vieta2} and related topics of
non-commutative geometry. Using the Gelfand--Retakh theory of
quasi-determinants, Bergvelt and Rabin in~\cite{rabin:duke} found an analog of
Cramer's formula in the supercase. The situation with Cramer's rule, i.e.,
calculating the inverse of a supermatrix, is a bit peculiar. At the first
glance one does not expect a role of the Berezinian similar to that of the
determinant in the classical case. However, this is true, though not so
straightforwardly (e.g., what should be the correct notion of a minor or an
adjunct? -- see in the main text). We give here a simple direct proof based on
the geometrical meaning of the Berezinian.

We would like to stress that  our methods throughout this paper are very
elementary.

The topics of our paper are intimately related with subtle questions concerning
rational and polynomial invariants of operators in superspaces. As it is known
(see below), the distinction between rational and polynomial invariants in the
$\ZZ$-graded situation is much sharper than in the classical case. On the other
hand, `rational' seems to be intrinsically related with `super'. For example,
every rational function $R(z)$ such that $R(0)=1$ can be viewed as the
characteristic function of a linear operator, $R(z)=R_A(z)$, its zeros and
poles corresponding to the bosonic and fermionic eigenvalues of $A$.  A pair of
polynomials $P,Q$ of degrees $p$ and $q$ can be viewed as the numerator and
denominator of such a characteristic function. One can show that their
resultant $\Res (P,Q)$ can be expressed via $\Tr\wed^kA$. It is the
(super)trace  of the representation corresponding to the $p\times q$
rectangular Young diagram $D$, $\Res (P,Q)=\Tr A_D$ (see in the main text).

Note that rational and polynomial invariants of supermatrices were first
considered in the pioneer works  on representations of Lie superalgebras by
Berezin (see references in ~\cite{berezin:antire}; some texts of 1975-77 were
incorporated into the English version of that posthumous book) and
Kac~\cite{kac:typical77}. They showed that all rational invariants of
supermatrices $p|q\times p|q$ can be expressed as rational functions of the
$p+q$ supertraces $\Tr A$, \ldots, $\Tr A^{p+q}$. In particular, one can
express polynomial invariants, though possibly non-polynomially. It was
discovered that not every polynomial on the diagonal matrices separately
symmetric in the `bosonic' and `fermionic' eigenvalues can be extended to a
polynomial invariant on matrices. In fact, it in general corresponds to a
rational invariant function  with a denominator of a special
appearance~\cite{berezin:superlaplace} (see the English version
of~\cite{berezin:antire}). In~\cite{berezin:superlaplace} Berezin  gave a
criterion for such function to be a polynomial, which was later clarified and
extended  to other Lie superalgebras by Sergeev~\cite{sergeev:bolg82} (see
also~\cite{sergeev:inv99}). There was an interesting sequel of works by Kantor
and Trishin~\cite{kantor-trishin:det, kantor-trishin:inv,
kantor-trishin:cayley}, in which the authors were concerned with clarifying the
relations in the (infinitely generated) algebra of polynomial invariants. In
particular, they found by a  method different from ours the
relations~\eqref{eqhankelc} for traces $\Tr \wed^k A$ for an arbitrary operator
in a $p|q$-dimensional space and came to analogs of the Cayley--Hamilton
identity. They did not consider expansions of rational functions. Their main
tool was analysis of  Young diagrams and the corresponding representations. The
coefficients of the expansion at infinity of the characteristic function (which
include the Berezinian) did not appear in these papers.

The recurrence relations linking $\wed^k V$ and $\Sigma^kV$ that we establish
in this paper, both as relations for traces and the relations in the
Grothendieck ring, are new. The explicit invariant formula for the Berezinian
as a ratio of two supertraces following from them, is also  new. Using these
relations more results can be obtained.  We hope that our approach allows to
reach better clarity of understanding of the Cayley--Hamilton theorem and
Cramer's rule in the supercase.

\subsection{Notation.} We use standard language of superalgebra and supergeometry.
Whenever it could not cause confusion, we drop the prefix `super', writing
`spaces', `traces', etc., instead of `super\-spaces',`super\-traces',
respectively.

\subsection{Acknowledgement.} We wish to thank A.~B.~Borovik, G.~Megy\-esi,
P.~N.~Pyatov, J.~Rabin, V.~S.~Retakh,  and A.~S.~Sorin for discussions at
various times, and Th.~Schmitt for sending us his paper~\cite{schmitt:ident}.
This work was reported at S.~P.~Novikov's seminar at the Steklov Institute in
Moscow. Special thanks go to him and the participants of the seminar for many
useful comments.

\section{Expansions of the Characteristic Function} \label{seccharfunct}

Let $A$ be an even linear operator acting in a finite-dimensional superspace
$V$. Denote $\dim V=p|q$. Consider the \textit{characteristic function} of this
operator,
   \begin{equation}\label{charfunction}
     R_A(z):=\Ber(1+zA)\,,
    \end{equation}
depending on a complex variable $z$. Here $\Ber$ denotes the Berezinian
(superdeterminant). If $M$ is an even invertible $p|q\times p|q$ supermatrix,
  $M=\left(\begin{matrix}
  M_{00} & M_{01} \\
  M_{10} & M_{11} \\
\end{matrix}\right)$, recall that
\begin{equation}\label{eqdefofber}
   \Ber M=\frac{\det \left(M_{00}-M_{01}M^{-1}_{11}M_{10}\right)}{\det
M_{11}}
\end{equation}
The Berezinian is a multiplicative function of matrices, hence it is
well-defined on linear operators.

Recall that for an even matrix (resp., operator), in the diagonal blocks
$M_{00}, M_{11}$ (resp., $A_{00}, A_{11}$) the matrix entries are even and in
the antidiagonal blocks $M_{01}, M_{10}$ (resp., $A_{01}, A_{10}$) the entries
are odd. In the sequel, when it cannot cause a confusion we do not distinguish
sharply operators and the corresponding matrices. Matrix elements can be viewed
either as belonging to a given $\ZZ$-graded (super)commutative ring or as free
generators. Classically this corresponds to considering an `individual' matrix
or a `general' matrix. Strictly speaking one should talk about `free modules'
over the ground ring instead of `vector spaces', but we shall not stress this
distinction.

Consider the expansion of the rational function $R_A(z)$ at zero:
\begin{equation}\label{eqexpansion}
R_A(z)=\sum_{k=0}^\infty c_k(A)z^k = 1+c_1 z+c_2z^2+\ldots \,.
\end{equation}

In the ordinary case (where the odd dimension  is equal to zero) the function $R_A(z)$
is a polynomial and the expansion~\eqref{eqexpansion} terminates. It is well known that
for a linear operator acting in a $p$-dimensional vector space $V$
\begin{equation*}\label{factformtextbook}
    \det (1+zA)=1+c_1z+\ldots+c_pz^p
\end{equation*}
where $c_k(A)=\Tr \wed^k A$ are the traces of the action of the operator $A$ in the
exterior powers $\wed^k V$. In particular, $c_1(A)=\Tr A$, $c_p(A)=\det A$. For $k>p$,
$c_k(A)=0$ as $\wed^k V=0$.

If the odd dimension of $V$ is not equal to zero, then $\Ber (1+zA)$ is no longer a
polynomial in $z$, but an analog of the formula above still holds:
\begin{prop}
There is an infinite power expansion
\begin{equation} \label{eqexpanzero}
    \Ber (1+zA)=\sum_{k=0}^\infty  c_k(A) z^k\quad \text{where $c_k(A)=
\Tr \wed^k A$}.
\end{equation}
\end{prop}

In~\eqref{eqexpanzero} $\wed^k A$ stands for the action of $A$ in the $k$-th exterior
power of the superspace $V$, where the exterior algebra $\wed (V)=\oplus \wed^k V$ is
defined as $T(V)/\langle v\otimes u+(-1)^{\vt\ut}u\otimes u\rangle$, $v,u$ being
elements of $V$. Parity in $\wed(V)$ (the $\ZZ$-grading) is naturally inherited from
$V$. There is no ``top'' power among $\wed^k V$, and the Taylor
expansion~\eqref{eqexpanzero} is infinite.

We denote the supertrace of a supermatrix  by the same symbol as the trace of an
ordinary  matrix. Recall that for an even supermatrix,
              $$
  \Tr M= \Tr\left(\begin{matrix}
  M_{00} & M_{01} \\
  M_{10} & M_{11} \\
 \end{matrix}\right)=\Tr M_{00}-\Tr M_{11}\,.
               $$

Expansion~\eqref{eqexpanzero} can be proved by considering diagonal matrices.
As far as we have managed to find out, this formula was  first obtained
in~\cite{schmitt:ident}.

The expansion of the characteristic function at infinity leads to traces of the wedge
products of the inverse matrix:
\begin{equation}\label{eqexpaninfinity}
    \Ber (1+zA)=\sum^{\infty}_{k=q-p} c_{-k}^*(A)z^{-k} \quad
    \text{where $c^*_{-k}(A)= \Ber A \cdot\Tr \wed^{p-q+k} A^{-1}$}.
\end{equation}
Formula~\eqref{eqexpaninfinity} follows from the equalities $\Ber (1+zA)=\Ber
A\,\Ber(A^{-1}+z)=z^{p-q}\Ber A\,\Ber(1+z^{-1}A^{-1})$ and
~\eqref{eqexpanzero}.  The geometric meaning of the
expansion~\eqref{eqexpaninfinity} is as follows. $\Ber A\cdot \Tr\wed^{p-k}
A^{-1}=\Tr \Sigma^k A$ is the trace of the representation of $A$ in the space
$\Sigma^k V :=\Ber V\otimes \wed^{p-k} V^* $. In the ordinary case, it would be
just a ``dual'' description of the same  $\wed^k V$; in the super case these
two spaces are essentially different.  Hence we get the following proposition.

\begin{prop}
There is an expansion at infinity
\begin{equation} \label{eqexpaninfinity2}
     \Ber (1+zA)=\sum_{k=q-p}^{\infty}  c_{-k}^*(A) z^{-k}\quad
     \text{where $c_{-k}^*(A)= \Tr \Sigma^{q-k}A$},
\end{equation}
which is a Taylor expansion when $p\leq q$ and a Laurent expansion when $p>q$. Here
$\Tr \Sigma^{q-k}A=\Ber A\cdot \Tr \wed^{p-q+k}A^{-1}$.
\end{prop}

Consider the coefficients $c_k(A)=\Tr\wed^k A$. They can be expressed as
polynomials via  $s_k(A)=\Tr A^k$. This follows from the Liouville formula
(hence basically from the multiplicativity of the Berezinian):
\begin{equation*}
\Ber (1+zA)=e^{\Tr\ln(1+zA)}=\exp \Bigl(z\Tr A- \frac{ z^2}{2} \,\Tr
A^2+\frac{z^3}{3} \,\Tr A^3+\ldots  \Bigr)\,.
\end{equation*}
Hence $c_k(A)$ can be expressed via $s_k(A)=\Tr A^k$ by the formulae
$c_k(A)=P_k(s_0(A),\dots,s_k(A))$, where $P_k$ are   classical Newton's
polynomials. For example, $c_0=s_0=1$,
          \begin{align*}
    c_1=s_1,\quad
    c_2=\frac{1}{2}\,(s_1^2-s_2),\quad
    c_3=\frac{1}{6}\,(s_1^3-3s_1s_2+2s_3),
           \end{align*}
etc., where $c_k=c_k(A)$, $s_k=s_k(A)$.
    There is a formula
\begin{equation}\label{newton}
 c_{k+1}=\frac{1}{k+1}\,(s_1 c_k- s_2 c_{k-1}+\ldots +(-1)^k
 s_{k+1}).
\end{equation}
These universal formulae linking $c_k(A)$  with $s_k(A)$ are true regardless whether
$V$ is a superspace or ordinary space.

For further considerations
it is convenient to define the following polynomials:
\begin{equation}\label{defofhampol}
   \H_k(z)=z^k-c_1z^{k-1}+c_2 z^{k-2}-\ldots+(-1)^kc_k,
\end{equation}
where   $k=0,1,2,\ldots$.  We shall refer to them as to the
\textit{Cayley--Hamilton polynomials}. (They appear with the relation to the
analog of the Cayley--Hamilton theorem which we discuss later. In the classical
case of an $n$-dimensional space,  $\pm\H_n(z)$ is the classical characteristic
polynomial $\det(A-z)$ if $c_k=c_k(A)$.) The following identities are
satisfied:
\begin{gather}\label{propertiesofham}
  \frac{dc_{k+1}(A)}{dA}=(-1)^k{\H}^A_k(A)\\
     \frac{1}{k+1}\Tr \bigl(A \H^A_k(A)\bigr)=(-1)^kc_{k+1}(A),
\end{gather}
Here   ${\H}^A_k(A)$ is the value of the polynomial~\eqref{defofhampol} where
$c_k=c_k(A)$ at $z=A$. The derivative $\dfrac{df(A)}{dA}$ of a scalar function
of a  matrix argument is defined as the matrix which satisfies
$\Bigl\langle\dfrac{df(A)}{dA}\,,
B\Bigr\rangle=\dfrac{d}{dt}f(A+tB)\big\vert_{t=0}$ for an arbitrary matrix $B$,
where the scalar product of matrices is given by $\langle A, B\rangle=\Tr AB$.
Formulae (\ref{newton}), (\ref{propertiesofham}) can be deduced by
differentiating the characteristic function $R_A(z)=\Ber(1+Az)$. Bearing in
mind that $d\Ber M=\Ber M\Tr (M^{-1}dM)$, we can come to the following
identities:
\begin{align*}
  \frac{d}{dz}\log R_A(z)&=\Tr \left(A\left(1+Az\right)^{-1}\right)=
    \sum_{k=0}^\infty (-1)^ks_{k+1}(A)z^k \,,\\
    \frac{d}{dA}\log R_A(z)&=\left(1+Az\right)^{-1}z=\sum_{k=0}^\infty
    (-1)^kz^{k+1}A^k\,.
\end{align*}
By writing $d\log R_A(z)$ as $(R_A(z))^{-1}dR_A(z)$ and comparing the power
series we arrive at (\ref{newton}), (\ref{propertiesofham}).

Unlike the polynomial functions $c_k(A)=\Tr\wed^k A$,  the coefficients
$c_k^*(A)=\Tr\Sigma^{q+k}A= \Ber A\cdot \Tr \wed^{p-q-k}A^{-1}$ are rational functions
of the matrix entries of $A$. In particular,
$$
c_{p-q}^*(A)=\Tr\Sigma^{p}A=\Ber A.
$$
Our task will be to give an expression for $c_k^*(A)$ in terms of polynomial
invariants of $A$.

\section{Recurrence Relations for  Traces of  Exterior
Powers}\label{secreltraces}

Recall that $\Sigma^k A$ denotes the representation of $A$ in the space
$\Sigma^k V=\Ber V\cdot \wed^{p-k}V^*$, thus $\Tr \Sigma^k A=\Ber A\cdot
\Tr\wed^{p-k} A^{-1}$.

By definition, $\Tr\Sigma^k A=0$ when $k>p\,\,$  and $\Tr \wed^k A=0$ when $k<0$.

In the purely even case ($q=0$, $\dim V=p$), the spaces $\wed^kV$ and $\Sigma^kV$ are
canonically isomorphic, $\Tr \wed^k A=\Tr \Sigma^k A$, and $c_p(A)=\det A$, $c_k(A)=0$
for $k>p$. We shall find out now what replaces these facts for a general
$p|q$-dimensional superspace.

Let us analyze the expansions of the characteristic function $R_A(z)$. One can see that
$R_A(z)$ is a fraction of the appearance
\begin{equation*}
R_A(z)=\frac{P(z)}{Q(z)}=\frac{1+a_1 z +a_2z^2+\ldots+ a_pz^p}{1+b_1 z +b_2z^2+\ldots+
b_qz^q}
\end{equation*}
where the numerator is a polynomial of degree $p$ and the denominator is a
polynomial of degree $q$. (Consider the diagonal matrices.) In principle the
degrees can be less than $p$ and $q$, and the fraction may be reducible.
However, for an operator ``in a general position'', this fraction is
irreducible and the top coefficients $a_p$, $b_q$ can be assumed to be
invertible. (A discussion of algebraic problems related with the notion of
``general position'' in this context can be found
in~\cite{kobayashi:nagamachi}. See also Section~\ref{seccayley}.) We shall use
the notation $\Rp_A(z)$ and $\Rm_A(z)$ for the numerator and denominator of the
fraction $R_A(z)$. Later we shall show how $\Rp_A(z)$ and $\Rm_A(z)$ can be
determined from the operator $A$.

From the well known connection between rational functions and recurrent
sequences (see Appendix), one can deduce the following facts:

\smallskip
(1) The coefficients $c_k(A)=\Tr \wed^k A$ of the expansion of $R_A(z)$ at
zero~\eqref{eqexpanzero} satisfy the recurrence relation of period $q$
\begin{equation}\label{eqrelationatzero}
b_0c_{k+q}+\ldots+b_qc_k=0
\end{equation}
for all $k>p-q$, where $b_0=1$. In particular, if $p<q$, then the
relation~\eqref{eqrelationatzero} holds for all $c_k$  including the zero values when
$p-q<k<0$.

(2) The coefficients $c_k^*(A)=\Tr \Sigma^{q+k}A$ of the expansion of $R_A(z)$  at
infinity~\eqref{eqexpaninfinity2} satisfy the same recurrence relation:
\begin{equation}\label{eqrelationatinfinity}
b_0c_{k}^*+\ldots+b_qc_{k-q}^*=0
\end{equation}
for all $k< 0$. In particular, if $p<q$, then the relation~\eqref{eqrelationatinfinity}
holds for all $c^*_k$ including the zero values when $p-q<k<0$.

(3) If $p<q$, then $c_k$ and $-c^*_k$ can be combined together into a single recurrent
sequence, for all $k\in \Z$:
\begin{equation}\label{eqchat}
\widehat{c}_k=
\begin{cases}
\phantom{-}c_k  & \text{if \quad $k\geq 0$}\\
\phantom{-}0 & \text{if \quad $p-q<k<0$} \\
-c^*_k  & \text{if \quad $k\leq p-q$}
\end{cases}
\end{equation}
The same holds in general: if one considers $c_k$ with sufficiently large positive $k$
and $-c^*_k$ with sufficiently large negative $k$, they fit into a single recurrent
sequence.

(4) Moreover, for arbitrary $p$ and $q$ the differences
\begin{equation*}
\dif_k=c_k-c^*_k.
\end{equation*}
satisfy the recurrence relation
\begin{equation}\label{eqdifrelation}
b_0\dif_{k+q}+\ldots+b_q\dif_k=0
\end{equation}
for \textit{all} values of $k\in\Z$ (notice that $c_k=0$ for $k<0$, $c^*_k=0$ for
$k>p-q$).

\smallskip
In particular, we have obtained the following fundamental theorem.

\begin{thm}\label{thmrecurrence}
 For an operator $A$ acting in $p|q$-dimensional vector space
 the differences
\begin{equation}\label{eqdiffer}
    \dif_k=c_k-c_k^*=\Tr\wed^k A-\Tr \Sigma^{q+k} A
\end{equation}
form a recurrent sequence with period $q$, for all $k\in \Z$. \qed
\end{thm}
In the classical case of $q=0$, all terms of the sequence~\eqref{eqdiffer} are zero and
$\Tr\wed^k A=\Tr \Sigma^{k}A$ or $\Tr\wed^k A=\det A\cdot \Tr \wed^{p-k} A^{-1}$ for
any operator $A$ which is a familiar equality.  In this case the spaces $\wed^k V$ and
$\Sigma^k V$ are canonically isomorphic. Theorem~\ref{thmrecurrence} actually suggests
a relation between spaces $\wed^k V$ and $\Sigma^{k+q} V$ for arbitrary $q$ (see
details in Section~\ref{secgroth}).

In~\eqref{eqdiffer} the terms $c_k=\Tr\wed^k A$ and $c_k^*=\Tr \Sigma^{q+k}A$
 can be both nonzero only in a finite range, for
$k=0,\ldots,p-q$ when $p>q$. Otherwise $\dif_k$ equals either $c_k(A)$ (for
$k\geq p-q+1$) or $-c_k^*(A)$ (for $k\leq -1$).  The relation \eqref{eqdiffer}
gives us a tool to express terms of the recurrent sequences
$c_k^*=\Tr\Sigma^{q+k} A$ and $c_k=\Tr\wed^k A$ via each other.

What actually happens, for large $k$, $\dif_k=c_k$, and they can be continued to the
left using~\eqref{eqdifrelation} to obtain $c_k^*(A)$, in particular $c_{p-q}^*(A)=\Ber
A$, as
\begin{equation*}
\Ber A =\Tr\wed^{p-q}A-\dif_{p-q}\,.
\end{equation*}
The  ``continuation to the left'' of $c_k(A)$ using the recurrence
relation~\eqref{eqrelationatzero} corresponds to the analytic continuation of the power
series~\eqref{eqexpanzero} representing the rational function $R_A(z)$ near zero.

\begin{exe}
If $p<q$, then   $\Tr\wed^k A$ and $-\Tr \Sigma^{q+k}A$ make a single recurrent
sequence for all $k$, so   $\dif_k=\widehat{c}_k$ in the notation above~\eqref{eqchat}.
Hence, in particular,
\begin{equation}
\Ber A=-\widehat{c}_{p-q}.
\end{equation}
\end{exe}
\noindent We give examples of calculations in the next section.

For linear recurrence relations with constant coefficients such as
~\eqref{eqrelationatzero} or~\eqref{eqdifrelation} it is possible to eliminate the
coefficients to obtain the relation ``in a closed form''. This is a standard method
based on the connection of recurrent sequences and rational functions with infinite
Hankel matrices (see, e.g.,~\cite{gantmacher:mat}). Recall that a Hankel matrix is one
with the entries $c_{ij}=c_{i+j}$. A recurrence relation for $c_k$ of period $q$
implies the vanishing of Hankel determinants of order $q+1$.

The statement \eqref{eqdiffer} of the Theorem can be reformulated in the
following way: the identity
\begin{equation}\label{eqhankelg}
\begin{vmatrix}
      \dif_k(A) & \dots & \dif_{k+q}(A) \\
      \dots & \dots & \dots \\
      \dif_{k+q}(A) & \dots & \dif_{k+2q}(A) \\
    \end{vmatrix}=0
\end{equation}
holds for all $k\in\Z$.

\begin{coro}\label{correcurc}
The identity
\begin{equation}\label{eqhankelc}
\begin{vmatrix}
      c_k(A) & \dots & c_{k+q}(A) \\
      \dots & \dots & \dots \\
      c_{k+q}(A) & \dots & c_{k+2q}(A) \\
    \end{vmatrix}=0
\end{equation}
holds for all $k> p-q$.
\end{coro}

\begin{Rem}
In works~\cite{berezin:superlaplace}, \cite{kac:typical77} appeared a system of
equations for $b_1,\ldots,b_q$  which is our equations~\eqref{eqrelationatzero}
(for $p\geq q$) with $k=p-q+1, \ldots, p$, but they did not consider recurrence
relations. The recurrence relations for $c_k=\Tr\wed^kA$, in particular the
identity~\eqref{eqhankelc}, appeared in~\cite{kantor-trishin:inv} and was then
interpreted in~\cite{kantor-trishin:cayley} by an analysis  of Young diagrams.
Compared to our work, in~\cite{kantor-trishin:inv} they came to the recurrence
relation for $c_k$  by pure combinatorics, using an explicit expression of
$c_k$  in terms of symmetric functions of the `bosonic' and `fermionic'
eigenvalues for a diagonal matrix, and not from  the characteristic function
$R_A(z)$, as we do here. Because of that, in the
works~\cite{kantor-trishin:inv, kantor-trishin:cayley} they never considered
the coefficients $c_k^*$; hence they could not see the general recurrence
relations involving both $c_k$ and $c_k^*$ that we establish here.
\end{Rem}

\section{Berezinian as a Rational Function of Traces}\label{secber}

As we established above, the coefficients $c_k(A)=\Tr \wed^k A$ for a linear
operator $A$ in a $p|q$-dimensional vector space $V$ satisfy relations
\eqref{eqrelationatzero} making them a $p|q$-recurrent sequence (see Appendix
for the necessary notions). Basing just on this fact we will give a recurrent
procedure for calculating the characteristic function $R_A(z)=\Ber (1+zA)$ and
the Berezinian of the operator $A$. Then we will present a closed formula for
$\Ber A$ using the relations \eqref{eqhankelg} of Theorem \ref{thmrecurrence}.

Let  $\c=\{c_n\}_{n\geq 0}$ be a  $p|q$-recurrent sequence such that $c_0=1$ . Denote
by $\R_{p|q}(z,\c)$ its generating function:
\begin{equation*}
\R_{p|q}(z,\c)=
\frac{1+a_1z+\ldots+a_pz^p}{1+b_1z+\ldots+b_qz^q}=1+c_1z+c_2z^2+\ldots\,.
\end{equation*}
The fraction $\R_{p|q}(z,\c)$ is defined by the first $p+q$ terms
$c_1,c_2,\dots,c_{p+q}$ of the sequence $\c$\,:
\begin{equation*}
\R_{p|q}(z,\c)= \R_{p|q}(z,c_1,\ldots,c_{p+q})\,.
\end{equation*}

In particular, if $A$ is a $p|q\times p|q$ matrix and $\{c_k\}$ is the sequence of the
traces of  exterior powers of the matrix $A$ ($c_k=c_k(A)=\Tr \wed^k A$),  then
$\R_{p|q}(z,\c)$ coincides with the characteristic function of  $A$:
\begin{equation}\label{asfunction}
    R_A(z)=\R_{p|q}(z,c_1(A),c_2(A),\ldots,c_{p+q}(A))\,.
\end{equation}

The rational functions $\R_{p|q}(z,\c)=\R_{p|q}(z,c_1,\dots,c_{p+q})$ have the
following properties:

\smallskip
(1) If $p\geq q$,  then the sequence $\c'$ defined by $c_k':=\frac{c_{k+1}}{c_1}$
(assuming that the coefficient $c_1$ is invertible) is a $p-1|q$-recurrent sequence and
\begin{gather}\label{firstrecurrentrule}
 \R_{p|q}(z,\c)=1+c_1z\,\R_{p-1|q}(z,\c'), \intertext{i.e.,}
\R_{p|q}(z,c_1,\ldots,c_{p+q})=1+c_1z\,
\R_{p-1|q}\left(z,\frac{c_2}{c_1},\ldots,\frac{c_{p+q}}{c_1}\right)\,.
\end{gather}

(2) The sequence $\c^\Pi=\{c_n^\Pi\}$  defined according to
\begin{equation*}
     1+c_1^\Pi z+c_2^\Pi z^2+ \ldots  =
     \frac{1}{1+ c_1z+c_2z^2+\ldots }\,,
\end{equation*}
for example
\begin{equation}\label{definitionoftilde}
     c_1^\Pi=-c_1, \  c_2^\Pi=-c_2+c_1^2, \  c_3^\Pi=-c_3+2c_1c_2-c_1^3,  \
     \,\ldots\,,
\end{equation}
is a $q|p$-recurrent sequence, and
\begin{equation}\label{secondrecurrentrule}
   \R_{p|q}\left(z,c_1,\ldots,c_{p+q}\right)=
   \frac{1}{\R_{q|p}(z,c_1^\Pi,\ldots,c_{p+q}^\Pi)}\,.
\end{equation}

(If  $A$ is a $ p|q\times p|q $ supermatrix and $A^\Pi$ is the parity reversed
$q|p\times q|p \,$ supermatrix, then  $c_k(A^\Pi)=c_k(A)^\Pi$.)

Using these properties one can express the rational function $\R_{p|q}$
corresponding to a $p|q$-recurrent  sequence via the rational function
$\R_{0|1}$ corresponding to a $0|1$-recurrent sequence, i.e., a geometric
progression. The steps are as follows. If $p<q$, we
apply~\eqref{secondrecurrentrule} to get a $p'|q'$-sequence with $p'>q'$. If
$p>q$, we repeatedly apply~\eqref{firstrecurrentrule} to decrease $p$.

\begin{ex}\label{excharfunctioncal}
Let $A$ be a $p|1\times p|1$ matrix. Then it follows from
\eqref{firstrecurrentrule} and \eqref{secondrecurrentrule} that
\begin{multline*}
    R_A(z)=\R_{p|1}\left(z,c_1(A),c_2(A),\ldots,c_{p+1}(A)\right)=\\
    1+c_1z\,\R_{p-1|1}\left(z,\frac{c_2}{c_1},\ldots,
    \frac{c_{p+1}}{c_1}\right)=\dots =\\
     1+c_1z+\dots+c_{p-1}z^{p-1}+
  c_pz^p\R_{0|1}\left(z,\frac{c_{p+1}}{c_p}\right)=\\
 1+c_1z+\dots+c_{p-1}z^{p-1}+\frac{c_pz^p}{1-\frac{c_{p+1}}{c_p}z}=\\
    1+c_1z+\dots+c_{p-1}z^{p-1}+\frac{c_p^2\,z^p}{{c_p}-{c_{p+1}}z}
\end{multline*}
\end{ex}

We can also deduce from here  formulae for the Berezinian. One can see from
\eqref{eqexpaninfinity} that for a $p|q\times p|q$ matrix $A$
\begin{equation}\label{berasfunctionofchar}
     \Ber A=\lim_{z\to \infty}z^{q-p}\,R_A(z)
\end{equation}
Let $\c=\{c_n\}$, $n\geq 0$, be  an arbitrary $p|q$-recurrent sequence such
that $c_0=1$ and let $\R(z,\c)$ be its generating function.   Then mimicking
\eqref{berasfunctionofchar} we define the \textit{Berezinian} of this sequence
by the formula
\begin{equation}\label{berofsequence}
   \B_{p|q}(\c)=\lim_{z\to \infty}z^{q-p}\,\R_{p|q}(z,\c).
\end{equation}
If $c_n=c_n(A)=\Tr\wed^k A$, then $\B_{p|q}(\c)=\Ber A$. From
\eqref{firstrecurrentrule} and \eqref{secondrecurrentrule} immediately follow
relations for $B_{p|q}$:
\begin{equation}\label{recforberezinians}
     \B_{p|q}(\c)=\B_{p|q}(c_1,\ldots,c_{p+q})=
        \begin{cases}
        c_1 \B_{p-1|q}\left(\c^\prime\right)\quad &\text{if $p\geq q+1$}\\
        1+c_1 \B_{p-1|q}\left(\c^\prime\right)\quad &\text{if $p=q$}\\
        \frac{1}{\B_{q|p}\left(\c^\Pi\right)}\quad &\text{if $p\leq q-1$}\\
        \end{cases}
\end{equation}
where the sequences $\c^\prime$ and $\c^\Pi$ are  defined as above.

Using these relations one can calculate the Berezinians of matrices in terms of traces.
Note that from these recurrent relations follows that if $p>q$ then for a
$p|q$-recurrent sequence $\c$,  its Berezinian $\B_{p|q}$ depends only on the
coefficients $c_{p-q},\ldots, c_{p}\,,\ldots,c_{p+q}$.

\begin{ex} \label{exber11} For a  $1|1\times 1|1$ matrix:
\begin{multline*}
    \Ber A=
\B_{1|1}(c_1(A),c_2(A))=1+c_1\B_{0|1}\Bigl(\frac{c_2}{c_1}\Bigr)=
        1+\frac{c_1}{\B_{1|0}\Bigl(\Bigl(\displaystyle
\frac{c_2}{{c_1}}\Bigr)^\Pi\Bigr)}=\\
1+\frac{c_1}{-\frac{c_2}{c_1}}=1-\frac{c_1^2}{c_2}=\frac{c_2-c_1^2}{c_2}=
         \frac{\Tr A^2+(\Tr A)^2}{\Tr A^2-(\Tr A)^2}
\end{multline*}
(we have applied Newton's formulae to get the last expression).
\end{ex}

\begin{ex}
\label{exberp1} For  a $p|1\times p|1$ matrix:
\begin{multline*}
    \Ber A=\B_{p|1}\left(c_{p-1}(A),\ldots, c_p(A),\ldots,c_{p+1}(A)\right)=
    c_{p-1}
\B_{1|1}\left(\frac{c_p}{c_{p-1}},\frac{c_{p+1}}{c_{p-1}}\right)\\
    =\frac{c_{p-1}c_{p+1}-c_p^2}{c_{p+1}}
\end{multline*}
\end{ex}

\begin{ex}\label{exber22}
For  a $2|2\times 2|2$ matrix:
\begin{multline*}
    \Ber A=\B_{2|2}\left(c_1(A),c_2(A),c_3(A),c_4(A)\right)=
    1+c_1\B_{1|2}\left(\frac{c_2}{c_1},\frac{c_3}{c_1},\frac{c_4}{c_1}\right)=\\
    =1+\frac{c_1}{\B_{2|1}\bigl(\bigl(
    \frac{c_2}{c_1}\bigr)^\Pi,\bigl(\frac{c_3}{c_1}\bigr)^\Pi,
     \bigl(\frac{c_4}{c_1}\bigr)^\Pi\bigr)}=\\
    1+\frac{c_1}{%
     \B_{2|1}\bigl(-\frac{c_2}{c_1},-\frac{c_3}{c_1}+\bigl(\frac{c_2}{c_1}\bigr)^2,
     -\frac{c_4}{c_1}+2\,\frac{c_2}{c_1}\frac{c_3}{c_1}-
     \bigl(\frac{c_2}{c_1}\bigr)^3\bigr)}=
\\ 1-\frac{c_1^2}{c_2\B_{1|1}\bigl(\frac{c_3}{c_2}-\frac{c_2}{c_1},
              \frac{c_4}{c_2}-\,\frac{2c_3}{c_1}+
              \bigl(\frac{c_2}{c_1}\bigr)^2\bigr)}=
1-\frac{c_1^2}{c_2\left(1-
\frac{\bigl(\frac{c_3}{c_2}-\frac{c_2}{c_1}\bigr)^2}{
\frac{c_4}{c_2}-\frac{2c_3}{ c_1}+\bigl(\frac{c_2}{c_1}\bigr)^2}\right)}
\end{multline*}
\end{ex}

The last expression can be further simplified, and in principle one can proceed
in this way to get the answer for arbitrary $q$, but at this point it is easier
to give a general formula. It  will reveal an unexpected link with classical
algebraic notions.

\section{Berezinian and Resultant}

Let $A$ be an even linear operator in a $p|q$-dimensional superspace. Consider
the relation \eqref{eqhankelg} of Theorem \ref{thmrecurrence} for $k=p-q$.
Recall that $\gamma_{p-q}=c_{p-q}-c^*_{p-q}$,  $\gamma_k=c_k$ for $k\geq p-q+1$
and $c^*_{p-q}=\Ber A$. Hence we have the following equalities:
\begin{multline*}
       0=\begin{vmatrix}
      \dif_{p-q}  & \dots & \dif_{p}  \\
      \dots & \dots & \dots \\
      \dif_{p}  & \dots & \dif_{p+q}  \\
    \end{vmatrix}=
    \begin{vmatrix}
      c_{p-q}  -\Ber A & \dots & c_{p}  \\
      \dots & \dots & \dots \\
      c_{p}  & \dots & c_{p+q}  \\
    \end{vmatrix}=\\
\begin{vmatrix}
      c_{p-q}  & \dots & c_{p}  \\
      \dots & \dots & \dots \\
      c_{p}  & \dots & c_{p+q}  \\
    \end{vmatrix}-
         \Ber A
        \begin{vmatrix}
      c_{p-q+2} & \dots & c_{p+1}  \\
      \dots & \dots & \dots \\
      c_{p+1}  & \dots & c_{p+q}  \\
    \end{vmatrix}
\end{multline*}
We arrive at the formula
\begin{equation}\label{eqberezinian1}
    \Ber A=\frac{\begin{vmatrix}
      c_{p-q} & \ldots & c_p \\
      \ldots & \ldots & \ldots \\
      c_p & \ldots & c_{p+q} \\
    \end{vmatrix}}{\begin{vmatrix}
      c_{p-q+2} & \ldots & c_{p+1} \\
      \ldots & \ldots & \ldots \\
      c_{p+1} & \ldots & c_{p+q} \\
    \end{vmatrix}}=\frac{|c_{p-q}\ldots c_p|_{q+1}}{|c_{p-q+2}\ldots
    c_{p+1}|_{q}}\,,
\end{equation}
where we used a short notation for Hankel determinants with subscripts denoting
their orders. Here as always $c_k=0$ for $k\leq -1$  and $c_0=1$.

Let us  make an important observation. By the Schur--Weyl character formula it
follows that the Hankel determinants appearing in the numerator and denominator
of formula \eqref{eqberezinian1} are nothing but the traces of the
representations of $A$ in the subspaces of tensors corresponding to certain
Young diagrams.

Indeed, denote by $D=D_{[\lambda_1,\dots,\lambda_s]}$ the Young diagram with
$s$ columns, such that the $i$-th column contains $\lambda_i$ cells,
$\lambda_1\geq \lambda_2\geq \dots\geq \lambda_s$. Let $V_D$ be an invariant
subspace in the tensor power $V^{\otimes\, N}$, $N=\lambda_1+\dots+\lambda_s$,
corresponding to the Young diagram $D=D_{[\lambda_1,\dots,\lambda_s]}$, and
$A_D$ be the representation of $A$ in $V_D$. Then the \textit{Schur--Weyl
formula} (see~\cite{weyl:classical}) tells that the trace of  $A_D$ is
expressed via the traces $c_k(A)=\Tr \wed^k A$ as the determinant of the
following $s\times s$ matrix:
\begin{gather*}\label{Shur}
    a_{ij}=c_{\lambda_i+j-i}(A)=
 \Tr \Lambda^{\lambda_i+j-i}A, \\
  \Tr   A_D=\det \left(a_{ij}\right).
\end{gather*}
It is known that the formula remains valid in the supercase (if trace means
supertrace). Let $D(r,s)$ be the rectangular Young diagram with $r$ rows and
$s$ columns. So $D(r,s)=D_{[\lambda_1,\dots,\lambda_s]}$ with $\lambda_i=r$ for
all $i$. One can see that for  $D=D(r,s)$ the `Schur determinant' $\Tr A_D$ is
equal to the Hankel determinant $|c_{r-s+1}\dots c_{r}|_{s}$ of order $s$, with
the inverted order of rows. In other words,  Hankel determinants appearing in
this paper can be interpreted as characters of tensor representations
corresponding to rectangular Young diagrams. Hence, in particular, our
formula~\eqref{eqberezinian1} for the Berezinian can be rewritten in the
following form
\begin{equation}\label{berformviaYoung}
\Ber A=(-1)^q\frac{\Tr A_{D(p,q+1)}}{\Tr A_{D(p+1,q)}}\,,
\end{equation}
the sign coming from the change of order of rows in the determinants.

\begin{Rem}
In the classical situation ($q=0$) when $c_k(A)$ are the elementary symmetric
functions of the eigenvalues of $A$, Schur's determinants corresponding to
Young diagrams (or partitions) when written as functions of these eigenvalues,
are special symmetric functions known as \textit{Schur functions}
(see~\cite{macdon:symm}); in the supercase the same Schur determinants when
expressed via the  eigenvalues are no longer classical symmetric Schur
functions but are combinations of functions that are separately symmetric in
the `bosonic' and `fermionic' eigenvalues. They should probably  be called
\textit{`super Schur functions'}.
\end{Rem}

\begin{ex} For a $2|3\times 2|3$ matrix we have
\begin{equation*}
\Ber A=\frac{\begin{vmatrix}
  0 & 1 & c_1&c_2 \\
  1 & c_1 & c_2 &c_3\\
  c_1 & c_2 & c_3 &c_4\\
  c_2 & c_3 & c_4 &c_5\\
\end{vmatrix}}{\begin{vmatrix}
  c_1 & c_2 &c_3 \\
  c_2 & c_3 &c_4\\
  c_3 & c_4 &c_5\\
\end{vmatrix}}=-
 \frac{\Tr A_{D(2,4)}}{\Tr A_{D(3,3)}}
\,.
\end{equation*}
\end{ex}

The formulae obtained above deserve to be called a theorem.

\begin{thm}
The Berezinian of a linear operator $A$ in a $p|q$-dimensional space is equal
to the ratio of the traces of the representations in the invariant subspaces of
tensors corresponding to the  rectangular Young diagrams $D(p,q+1)$ and
$D(p+1,q)$
\begin{equation}\label{eqberezinian2}
    \Ber A= \frac{|\Tr\wed^{p-q}A\,\ldots\, \Tr\wed^pA|_{q+1}}%
    {|\Tr\wed^{p-q+2}A\,\ldots\,
    \Tr\wed^{p+1}A|_{q}}=\pm\,\frac{\Tr A_{D(p,q+1)}}{\Tr A_{D(p+1,q)}}
 \,.
\end{equation}
Here at the right hand side stand the Hankel determinants of orders $q+1$ and
$q$ made of the traces of  exterior powers of the operator $A$. \qed
\end{thm}

What is the meaning  --- as polynomial invariants of $A$ --- of the
determinants $\Tr A_{D(p,q+1)}$ and $\Tr A_{D(p+1,q)}$  appearing as the
numerator and denominator in formula~\eqref{eqberezinian2}?

\begin{de} Define the following functions of $A$:
\begin{align}
   \Berp A&:=\l_1\ldots\l_p \prod_{i,\a} (\l_i-\m_{\a}) \,,\\
    \Berm A&:=\m_1\ldots\m_q \prod_{i,\a} (\l_i-\m_{\a}).
\end{align}
\end{de}

We assume for a moment that $A$ can be diagonalized and  $\l_i$, $\m_{\a}$,
$i=1,\ldots,p$, $\a=1,\ldots,q$ stand for its eigenvalues. So
$$\Ber
A=\dfrac{\l_1\ldots\l_p}{\m_1\ldots\m_q}=\dfrac{\Berp A}{\Berm A}\,.
$$
We shall immediately see that $\Ber^\pm A$   make sense for all $A$.

Denote the product $\prod_{i,\a} (\l_i-\m_{\a})$ by $R$ or $R(A)$.  If
$\Rp_A(z)$ and $\Rm_A(z)$ stand for the numerator and denominator of the
characteristic function $R_A(z)$, then it is easy to check that $R$ is the
classical Silvester's resultant for the polynomials $\Rp_A(z)$ and $\Rm_A(z)$,
$R=\Res(\Rm_A(z), \Rp_A(z))$.

\begin{prop}\label{propresultant}
The resultant of $\Rp_A(z)$ and $\Rm_A(z)$ can be expressed by the following
formula:
\begin{multline}\label{eqresrprm}
    R=\Res(\Rm_A(z), \Rp_A(z))=\prod_{i,\a} (\l_i-\m_{\a})=\\(-1)^{q(q-1)/2}|c_{p-q+1}\ldots
    c_{p}|_q=\Tr A_{D(p,q)}.
\end{multline}
\end{prop}
\begin{proof}
The Hankel determinant in the \RHS of \eqref{eqresrprm} vanishes when
$\l_i=\m_{\a}$ for any pair $i,\a$. This follows from our recurrence
relation~\eqref{eqhankelc} applied a $(p-1|q-1)$-dimensional space. Hence
$|c_{p-q+1}\ldots c_{p}|_q$ is divisible by the resultant. As polynomials in
$\l_i$, $\m_{\a}$ they have the same degree $pq$, hence they must coincide up
to a numerical factor, which can be checked, for example, by setting all
$\m_{\a}=0$.
\end{proof}

It follows that $R=R(A)$ is a polynomial in the matrix entries of $A$.

Note that the statement of Proposition~\ref{propresultant} is present in
Berezin's paper~\cite{berezin:superlaplace}.

\begin{thm}
The following equalities hold:
\begin{align}
   \Berp A&=\l_1\ldots\l_p \prod_{i,\a} (\l_i-\m_{\a}) = |c_{p-q}\ldots c_p|_{q+1} \label{eqberphank}\\
    \Berm A&=\m_1\ldots\m_q \prod_{i,\a} (\l_i-\m_{\a}) =|c_{p-q+2}\ldots
    c_{p+1}|_{q}\,, \label{eqbermhank}
\end{align}
i.e., $\Berp A$ and $\Berm A$ give exactly the top and bottom of the expression
for $\Ber A$ in formula~\eqref{eqberezinian2}.
\end{thm}

\begin{proof}
Indeed,  $\l_1\ldots\l_p$ and $\m_1\ldots\m_q$ are equal, respectively, to the
coefficients $a_p$ and $b_q$ in $\Rp_A(z)$ and $\Rm_A(z)$. In general, all the
coefficients $a_i$, $b_k$ can be obtained from $c_k$, $k=1,\ldots, p+q$, by
solving simultaneous equations, with the determinant of the system being
exactly $R$. Therefore, all coefficients $a_i$, $b_k$ have the appearance of a
polynomial in $c_k$ divided by the same denominator $R=\pm|c_{p-q+1}\ldots
c_{p}|_q=\Tr A_{D(p,q)}$. Formulae~\eqref{eqberphank} and \eqref{eqbermhank}
follow by a direct application of Cramer's rule. (In particular, this yields
another proof of the expression for the Berezinian~\eqref{eqberezinian2}.)
\end{proof}

From the proof, in particular, follows that the polynomials $\Rp(z)$ and
$\Rm(z)$ are  defined if the resultant $R=|c_{p-q+1}\ldots c_{p}|_q$ is
invertible.

Notice that the top and bottom of  the standard definition of the Berezinian
given by fraction~\eqref{eqdefofber} are non-invariant and non-polynomial
functions of the matrix; the products $\l_1\ldots\l_p$ and $\m_1\ldots\m_q$ are
invariant, but non-polynomial (and defined not explicitly as functions of the
matrix entries). The functions $\Ber^{\pm} A$ are polynomial invariants, and,
as one can see, they are the ``minimally possible'' modifications of the
products of eigenvalues with this property.

We have four remarkable Hankel (or Schur) determinants in this paper: $\Tr
A_{D(p,q)}$, $\Tr A_{D(p+1,q)}$, $\Tr A_{D(p,q+1)}$ and $\Tr A_{D(p+1,q+1)}$;
the first  being the resultant $R$, the last  giving the
identity~\eqref{eqhankelc} of the smallest degree, and the two in the middle
arising in the formula for the Berezinian~\eqref{eqberezinian2}.

\begin{Rem}
As a by-product of Proposition~\ref{propresultant} we have the following
formula for the resultant of two polynomials:
\begin{equation}\label{eqresultant}
    \Res(Q,P)=  \begin{vmatrix}
      c_{p-q+1} & \ldots & c_p \\
      \ldots & \ldots & \ldots \\
      c_p & \ldots & c_{p+q-1} \\
    \end{vmatrix}
\end{equation}
where $P(z)=a_pz^p+\ldots +1$, $Q(z)=b_qz^q+\ldots +1$, and the coefficients
$c_k=c_k(Q,P)$ are defined as follows:
\begin{equation}
    c_k(Q,P)=\sum_{i+j=k}a_i\tau_j(-1)^j
\end{equation}
where $\tau_j$ are the complete symmetric functions of the roots of $Q$. The
\RHS{} of \eqref{eqresultant} can be interpreted as  the (super)trace $\pm\Tr
A_D(p,q)$, where $A$ is an operator in a $p|q$-dimensional space associated
with the pair of polynomials $P,Q$, so that $R_A=\dfrac{P}{Q}$.
\end{Rem}

\section {Rational and Polynomial Invariants  and the Cayley--Hamilton
Identity}\label{seccayley}

\def\var#1#2 {#1_1,\dots,#1_#2 }

In the previous section we obtained  explicit formulae expressing the
Berezinian of a linear operator $A$ as rational function of traces. The
Berezinian is an example of a rational invariant function on supermatrices. Let
us briefly review general facts concerning such functions. This will be applied
to the analysis of the analog  of the Cayley--Hamilton theorem.

In the classical case invariant rational functions  $F(A)$ on $p\times p$
matrices,  $F(A)=F(C^{-1}AC)$, are in a $1-1$ correspondence with rational
symmetric functions $f(\var\l p)$ of $p$ variables, the eigenvalues of $A$. The
same is true for polynomial functions, due to the fundamental theorem on
symmetric functions and to the fact that the elementary symmetric polynomials
$\s_k(\l)$ (or the power sums $s_k(\l)$) are restrictions of the polynomial
functions of matrices $\Tr\wed^k A$ (resp., $\Tr A^k$).

This is not the case for $p|q\times p|q$ matrices, where arises a sharp
distinction between rational and polynomial invariants.

Every invariant rational function $F(A)$ on $p|q\times p|q$ matrices, i.e.,
$F(A)=F(C^{-1}AC)$ for every even invertible matrix $C$, defines a function
$f(\l_1,\ldots,\l_p,\m_1,\ldots,\m_q)$  of the eigenvalues of  $A$, with $\l_i$
corresponding to even eigenvectors and  $\m_{\a}$   to  odd eigenvectors,
symmetric separately in the variables $\var\l p$ and $\var\mu q$ (because even
and odd eigenvectors  cannot be permuted by a similarity transformation).

\begin{prop}\label{everythexpressedvia}
Every rational $S_p\times S_q$-invariant function of $\l_i,\m_{\a}$ can be
expressed as a rational function of  the polynomials $c_1,\ldots, c_{p+q}$ or
$s_1,\ldots, s_{p+q}$, where $c_k(\l,\m)=\Tr \wed^kA$, $s_k(\l,\m)=\Tr A^k$.
{\rm (Traces are supertraces).}
\end{prop}

\begin{ex}\label{strange}
Consider  the $S_1\times S_1$-invariant polynomial $f(\l,\mu)=\l+\mu$.  We have
\begin{align}\label{simplecountrexample}
        \lambda+\mu=
        \frac{\lambda^2-\mu^2}{\lambda-\mu}=\frac{s_2}{s_1}=\frac{c_1^2-c_2}{c_1}\,,
\end{align}
therefore it corresponds to a rational invariant function on $1|1\times 1|1$
matrices.
\end{ex}

We see that $S_p\times S_q$-invariant polynomials do not necessarily extend to
invariant polynomials of matrices.

Proposition~\ref{everythexpressedvia} (Berezin~\cite{berezin:superlaplace},
\cite[p. 315]{berezin:antire}, Kac~\cite{kac:typical77}) immediately follows
from considerations of the previous section, as all $S_p\times S_q$-invariant
functions of $\l_i,\m_{\a}$ are expressed via the elementary symmetric
functions of $\l_i$ and $\m_{\a}$, i.e., the coefficients $a_k$, $b_k$ of the
numerator and denominator of the characteristic function $R_A(z)$, which are
rational functions of $c_1,\ldots,c_{p+q}$. Moreover, \textit{for $S_p\times
S_q$-invariant polynomials $f(\l,\m)$ it follows that the corresponding
rational invariant functions $F(A)$ can be written as fractions with the
numerator being a polynomial invariant function of $A$ and the denominator
being a power of the resultant $R=R(A)$}.

The following non-trivial statement holds.

\begin{prop}[Berezin, Sergeev] For a $S_p\times S_q$-invariant polynomial $f(\l,\m)$
three conditions are equivalent: {\rm (a)} the equation
\begin{equation}\label{serg}
   \left(\frac{\p f}{\p \lambda_i}+\frac{\p f}{\p \mu_j}\right)\big\vert_{\lambda_i=\mu_j}=0,
\end{equation}
is satisfied; {\rm (b)} $f(\l,\m)$ extends to a polynomial invariant on
matrices; {\rm (c)} $f(\l,\m)$ can be expressed as a polynomial of a finite
number of functions $c_k(\l,\m)$, $k=0,1,2,3,\dots$ {\rm (or $s_k(\l,\m)$,
$k=0,1,2,3,\dots$)}.
\end{prop}
The implication (c)$\Rightarrow $(b) is obvious, the implication
(b)$\Rightarrow $(a) can be deduced from the invariance condition, the
implication (a)$\Rightarrow $(c) is the most technical part.
(See~\cite{berezin:superlaplace}, \cite[p. 294]{berezin:antire},
\cite{sergeev:bolg82}, \cite{sergeev:inv99}.)

\begin {ex}
\label{expolyratio} The $S_1\times S_1$-invariant polynomial
$f(\lambda,\mu)=\mu^N(\lambda-\mu)$ satisfies~\eqref{serg} and is in fact equal
to the polynomial $(-1)^Nc_{N+1}(A)$. It cannot be expressed as a polynomial in
$c_1,\dots,c_k$ if $k\leq N$. On the other hand, in full accordance with
Proposition~\ref{everythexpressedvia}, we can express it rationally via
$c_1,c_2$:
\begin{equation*}
\mu^N(\lambda-\mu)=(-1)^Nc_{N+1}(A)=
   \frac{c_2^N}{c_1^{N-1}}.
\end{equation*}
\end{ex}

Example~\ref{expolyratio} demonstrates that, differently from the classical
case, the algebra of polynomial invariants on supermatrices is not finitely
generated (no a priori number of $c_k$ is sufficient) and is not free (the
generators $c_k$, $k=1,2, \ldots\,$ satisfy an infinite number of
relations~\eqref{eqhankelc}).

\begin{Rem}
It would be interesting to describe the class of  invariant rational functions
on $\l_i,\m_{\a}$ that obey equation~\eqref{serg}. For example, the
characteristic function $R_A(z)$ and the Berezinian  $\Ber A$ belong to this
class. Hence it contains  products of polynomial invariants with arbitrary
powers of the Berezinian.
\end{Rem}

Now let us turn to the  Cayley--Hamilton theorem.

For an operator $A$ in a $p|q$-dimensional space it is clear that it
annihilates the polynomial $\HC_A(z)=\prod(\l_i-z)(\m_{\a}-z)$, where
$\l_i,\m_{\a}$ stand for the eigenvalues of $A$ as above, and one can see that
every polynomial annihilating a  generic operator $A$ is divisible by
$\HC_A(z)$, exactly as it is in the classical case. Hence, the polynomial
$\HC_A(z)$ is a minimal polynomial for generic operators. `Generic' means here
that all the differences of the eigenvalues, $\l_i-\l_j$, $\l_i-\m_{\a}$,
$\m_{\a}-\m_{\b}$, are invertible. In particular, $R=\Res(\Rm_A,\Rp_A)$ is
invertible and $\R^\pm_A(z)$ make sense. This \textit{`classical characteristic
polynomial'} or \textit{`Cayley--Hamilton polynomial'} of $A$, is expressed in
terms of the characteristic function $R_A(z)$ as
\begin{multline}\label{hameq1}
    \HC_A(z)=(-z)^{p+q}\Rp_A\Bigl(-\frac{1}{z}\Bigr)\Rm_A\Bigl(-\frac{1}{z}\Bigr)=\\
    \bigl(a_p-a_{p-1}z+\ldots+(-1)^pz^p\bigr)\bigl(b_q-b_{q-1}z+\ldots+(-1)^qz^q\bigr).
\end{multline}
Since the coefficients of $\R_A^\pm(z)$ are rational invariant functions of
$A$, with the denominator $R=\Res(\Rm_A,\Rp_A)=\Tr A_{D(p,q)}$, it follows that
the coefficients of $\HC_A(z)$, too, are  rational (not polynomial) invariant
functions of $A$, with  denominators $R$ or $R^2$.

\begin{ex}\label{naiveexample}
Consider a linear operator $A$ in a $p|1$-dimensional vector space $V$. Let us
calculate for it the  polynomial $\HC_A(z)$, which is here
$\HC_A(z)=(\lambda_1-z)\ldots(\lambda_p-z)(\mu-z)$. From
Example~\ref{excharfunctioncal} we get
\begin{multline*}
   R_A(z)=1+c_1z+\ldots+c_{p-1}z^{p-1}+
   \dfrac{c_p}{1-\dfrac{c_{p+1}}{c_p}z}\,z^p=\\
   \left(1+\frac{c_1c_p-c_{p+1}}{c_p}\,z +
   \frac{c_2c_p-c_1c_{p+1}}{c_p}\,z^{2}+
   \ldots+
   \frac{c_{p}c_p-c_{p-1}c_{p+1}}{c_p}\right)\times \\
   \left(1-\frac{c_{p+1}}{c_p}\,z\right)^{\!\!-1}
\end{multline*}
where $c_k=c_k(A)=\Tr \wed^k A$. Hence
\begin{multline*}
\HC_A(z)=(-1)^{p+1}\left(z^p-\frac{c_1c_p-c_{p+1}}{c_p}\,z^{p-1}+
\frac{c_2c_p-c_1c_{p+1}}{c_p}\,z^{p-2}-\right.\\
\left.\ldots+(-1)^p\,
\frac{c_{p}c_p-c_{p-1}c_{p+1}}{c_p}\right)\left(z+\frac{c_{p+1}}{c_p}\right)
\end{multline*}
and after simplification using the identity $c_pc_{p+2}-c_{p+1}^2=0$ we get
\begin{equation}\label{eqnaivehc}
    \HC_A(z)=\sum_{k=0}^{p+1}
    (-1)^{p+1-k}\,\frac{c_kc_p-2c_{k-1}c_{p+1}+c_{k-2}c_{p+2}}{c_p}\,z^{p+1-k}
\end{equation}
where as always $c_k=0$ for $k<0$. Notice that here $R=c_p$, and it appears in
the denominator in the final answer in the first power, not as $R^2$ as one
might expect, due to identities for $c_k$. We will see that this is the general
case.
\end{ex}

By multiplying $\HC_A(z)$ by its denominator we can get an annihilating
polynomial with the coefficients which are polynomial invariant functions of
the matrix entries of $A$. The advantage of such a polynomial is that it will
be an annihilating polynomial for arbitrary operators, not necessarily generic.
Notice that a minimal polynomial for generic operators is unique up to a factor
$R^N$.

\begin{ex} (Example \ref{naiveexample} continued.)
Multiplying both sides of~\eqref{eqnaivehc} by $c_p$ we obtain the polynomial
\begin{equation}\label{eqnaivehcpol}
    \tilde\HC_A(z)=\sum_{k=0}^{p+1}
    (-1)^{p+1-k} \left(c_kc_p-2c_{k-1}c_{p+1}+c_{k-2}c_{p+2}\right) z^{p+1-k}\,,
\end{equation}
which annihilates an arbitrary operator $A$ in a $p|1$-dimensional space and
whose coefficients are polynomial invariants of $A$.
\end{ex}

Let us show that the `naive'  characteristic  polynomial discussed above
follows also from the recurrence relations of Theorem~\ref{thmrecurrence}. A
method of constructing a `Cayley--Hamilton identity' from a relation on traces
was given in~\cite{kantor-trishin:cayley}. Below we shall use that method and
then show that the final answer can be identified with the naive formula
\eqref{hameq1} up to a factor.

If $A$ is an even linear operator in a $p|q$-dimensional vector space, then, in
particular, the traces of its exterior powers obey relations \eqref{eqhankelc}
for all $k>p-q$. For $k=p-q+1$ we have
\begin{equation}\label{eqhankelpol2}
\begin{vmatrix}
      c_{p-q+1}(A) & \dots & c_{p+1}(A) \\
      \dots & \dots & \dots \\
      c_{p+1}(A) & \dots & c_{p+q+1}(A) \\
    \end{vmatrix}=|c_{p-q+1}(A)\dots c_{p+1}(A)|_{q+1}=0\,.
\end{equation}
This is a scalar equation valid for any even matrix in a $p|q$-dimensional
space. Hence, by differentiating it one obtains a matrix identity (compare with
a formal differential calculus developed in ~\cite{kantor-trishin:cayley}).

In the classical case when $A$ is a linear operator in a $p$-dimensional vector
space ($q=0$), the relation~\eqref{eqhankelpol2} reduces to $c_{p+1}(A)\equiv
0$. Differentiating this identity gives exactly the vanishing of the
Cayley--Hamilton polynomial $\P_p(z)$ with  $c_k=c_k(A)$ at $z=A$, i.e., the
classical Cayley--Hamilton theorem.

For arbitrary $q$, by taking the derivative of~\eqref{eqhankelpol2} and
applying~\eqref{propertiesofham}, we get the equality
\begin{equation}\label{eqhcidentity2}
    \sum_{r=p-q+1}^{p+q+1} (-1)^{r-1} F^A_{r}\H^A_{r-1}(A)=0
\end{equation}
where we denote by $F_r$ the partial derivative of the Hankel determinant
$|c_{p-q+1}\dots c_{p+1}|_{q+1}$,
\begin{equation}\label{superHC}
   F_r =\dfrac{\p}{\p c_r}\, |c_{p-q+1} \dots c_{p+1}|_{q+1}\,,
\end{equation}
and by $F^A_r$ its value when $c_k=c_k(A)$.  Define a polynomial in $z$ of
degree $p+q$, with coefficients polynomially depending on $c_k$:
\begin{equation}\label{superham2}
          \tilde\HC(z):=    \sum_{r=p-q}^{p+q}(-1)^r F_{r+1}{\H}_r(z)\,.
\end{equation}
We shall write $\tilde\HC(z)=\tilde\HC_A(z)$ if $c_k=c_k(A)$. It follows that
$\tilde\HC_A(z)$ is an annihilating polynomial for $A$.

\begin{ex}
Let us make a calculation for $p|1\times p|1$ matrices. We have the identity
\begin{equation}\label{cancels}
     c_{p}c_{p+2}-c_{p+1}^2\equiv 0.
     \end{equation}
By differentiating we get $F_p=c_{p+2}$, $F_{p+1}=-2c_{p+1}$,  $F_{p+2}=c_p$.
Thus
$$
\tilde\HC_A(z)=(-1)^{p-1}F_{p}\H_{p-1}(z)+(-1)^{p}F_{p+1}\H_{p}(z)+(-1)^{p+1}F_{p+2}\H_{p+1}(z).
$$
After substituting the expressions~\eqref{defofhampol} for $\H_r(z)$ and
collecting similar terms we immediately get
\begin{equation}
    \tilde\HC_A(z)=\sum_{k=0}^{p+1}(-1)^{p+1-k}\left(c_kc_p-2c_{k-1}c_{p+1}+c_{k-2}c_{p+2}\right)
    z^{p+1-k}
\end{equation}
which precisely coincides with $R\cdot \HC_A(z)$ of Example~\ref{naiveexample}.
\end{ex}

Now we shall prove in general that by differentiating the identity for
traces~\eqref{eqhankelpol2}  one  arrives at a multiple of the `classical'
characteristic polynomial $\HC_A(z)$. Indeed, for generic matrices, $\HC_A(z)$
is a minimal polynomial, and any annihilating polynomial for $A$ is divisible
by $\HC_A(z)$. Consider the polynomial $\tilde\HC_A(z)$ defined in
~\eqref{superham2}. Dividing it by $\HC_A(z)$ we get $\tilde\HC_A(z)=c\cdot
\HC_A(z)$, where $c$ is a constant (as both polynomials are of the same
degree). To calculate $c$ compare the top coefficient in $\tilde\HC_A(z)$,
which is $(-1)^{p+q}F_{p+q+1}$, with that of $\HC_A(z)$, which is $(-1)^{p+q}$.
We have directly
\begin{equation*}
    F_{p+q+1} =\dfrac{\p}{\p c_{p+q+1}}\,
    \begin{vmatrix}
      c_{p-q+1} & \ldots & c_{p+1} \\
      \ldots & \ldots & \ldots \\
      c_{p+1} & \ldots & c_{p+q+1} \\
    \end{vmatrix}=
    \begin{vmatrix}
      c_{p-q+1} & \ldots & c_{p} \\
      \ldots & \ldots & \ldots \\
      c_{p} & \ldots & c_{p+q-1} \\
    \end{vmatrix}=R\,.
\end{equation*}
It follows that $c=R$. (We see that remarkably, $R$, not $R^2$, is the common
denominator of the fractions that are the coefficients of $\HC_A(z)$.) We
arrive at the following proposition.
\begin{prop}
The polynomial $\tilde\HC_A(z)$ defined by formula~\eqref{superham2} where
$c_k=c_k(A)$, is an anni\-hi\-lating polynomial for any operator $A$ in a
$p|q$-dimen\-sion\-al space. Its coefficients are invariant polynomial
functions of $A$. For generic operators, $\tilde \HC_A(z)$ is a minimal
polynomial, which divides all annihilating polynomials for $A$. The identity
holds:
\begin{equation}
    \tilde\HC_A(z)=R\cdot \HC_A(z)\,,
\end{equation}
where $\HC_A(z)=\prod (\l_{i}-z)(\m_{\a}-z)$ is the naive characteristic
polynomial, with rational coefficients, and
$R=\prod(\l_i-\m_{\a})=\Res(\Rm_A,\Rp_A)$.
\end{prop}

One can call the polynomial $\tilde \HC_A(z)$, with polynomial coefficients, a
\textit{`modified characteristic polynomial'}. In the classical situation,
holds $\HC_A(z)=\tilde\HC_A(z)=\det(A-z)$.

\section{Recurrence Relations in the Grothendieck Ring}\label{secgroth}

Recurrence relations for the traces of   exterior powers of an operator $A$ in
a $p|q$-dimensional superspace hold good for any operator,  their form being
independent of the operator. Such universal relations for   traces  suggest the
existence of underlying relations for the spaces themselves such as in the case
of $q=0$ the equality $\wed^kV=0$  when $k>p$. We shall deduce these relations
now.

First of all, let us explain in which sense we may speak about recurrence
relations for vector spaces. They hold in a suitable Grothendieck ring. One can
consider the Grothendieck ring of the category of all finite-dimensional vector
superspaces (i.e., $\ZZ$-graded vector spaces). This ring is isomorphic to
$\Z[\Pi]/\langle\Pi^2-1\rangle$, which is the ring where dimensions of
superspaces take values. An equality in this ring means just the equality of
dimensions. Alternatively, one can fix a superspace $V$ and consider the
Grothendieck ring of the category of all finite-dimensional superspaces with an
action of the supergroup $\GL(V)$, i.e., the Grothendieck ring of the
finite-dimensional representations of $\GL(V)$. Equality of two ``natural''
vector spaces like spaces of tensors over $V$ in this ring should mean the
existence of an isomorphism  commuting with the action of $\GL(V)$.

As a starting point we use the following relation, which holds for any superspace $V$:
\begin{subequations}\label{eqgroth1}
\begin{align}\label{eqgroth1a}
\wed_z(V)\cdot S_{-z}(V)&=1 \,,\\
\intertext{which one might prefer to rewrite as} \wed_z(V)\cdot \wed_{-z\Pi}(\Pi V)&=1
\label{eqgroth1b}
\end{align}
(for a proof it is sufficient to consider one-dimensional spaces).
\end{subequations}
Here $\wed_z(V)=\sum z^k\wed^k V =1+zV+z^2\wed^2V+\ldots\,$, etc. These are
power series in either of the Grothendieck rings described above. We  denote
the class of a vector space the same as the space itself. Notice that the unity
$1$ is the class of the main field. Equalities~\eqref{eqgroth1} hold in both
senses. For example, expanding in $z$ one gets $V-V=0$, $S^2V+\wed^2V-V\otimes
V=0$, etc.

Now, for a superspace $V$ we have $V=V_0\oplus V_1$ where $V_0$ is purely even and
$V_1$ is purely odd. We can rewrite this as $V=U\oplus \Pi W$ where both $U$, $W$ are
purely even vector spaces. It follows that $\wz(V)=\wz(U)\wz(\Pi W)$, therefore
by~\eqref{eqgroth1b}
\begin{equation}\label{eqgrothfrac}
    \wz(V)=\frac{\wz(U)}{\wzp(W)}=
    \frac{1+zU+z^2\wed^2U+\ldots+z^p\wed^pU}{1-z\Pi
    W+z^2\wed^2W-\ldots+(-z)^q\Pi^q\wed^qW}\,.
\end{equation}
Note that though $U$ and $W$ with their exterior powers do not belong to the
ring of representations of $\GL(V)$, they can be thought of as ideal elements
that can be adjoined to it, or, which is the same, as elements of the
representation ring of the block-diagonal subgroup $\GL(U)\times \GL(W)\subset
\GL(V)$. We see that the power series $\wz(V)$  represents a rational function
with the numerator of degree $p$ and denominator of degree $q$. Denote it by
$R_V(z)$; it replaces  the characteristic function $R_A(z)=\Ber(1+zA)$ of our
previous analysis. $R_A(z)$ can be viewed as the character of $R_V(z)$, for the
ring of representations of $\GL(V)$.

We can apply to $R_V(z)$ the same reasoning as to $R_A(z)$ above and conclude that the
exterior powers $\wed^kV$ for a $p|q$-dimensional vector space $V$ satisfy a recurrence
relation of period $q$
\begin{equation}\label{eqgrothrelation1}
b_0\wed^{k+q}V+\ldots + b_q\wed^kV=0
\end{equation}
for all $k\geq p-q+1$. Here $b_i=(-\Pi)^i\wed^iW$. Evidently in the classical case of
$q=0$ this reduces to $\wed^kV=0$ for $k\geq p+1$. The  relations for $c_k(A)=\Tr\wed^k
A$ then follow from~\eqref{eqgrothrelation1}.

As in Section~\ref{secreltraces}, it is possible to eliminate the coefficients
$b_i=(-\Pi)^i\wed^iW$ from the recurrence relations~\eqref{eqgrothrelation1} and
express them in a closed form using Hankel determinants. We arrive at the following
theorem.

\begin{thm}
\label{thmhankelgroth} For an arbitrary $p|q$-dimensional vector space $V$ the
following Hankel determinants vanish:
\begin{equation}\label{eqgrothrelation2}
    \begin{vmatrix}
      \wed^kV & \dots & \wed^{k+q}V \\
      \dots & \dots & \dots \\
      \wed^{k+q}V & \dots & \wed^{k+2q}V \\
    \end{vmatrix}=0
\end{equation}
for all $k\geq p-q+1$. \qed
\end{thm}

Notice that the expression of the recurrence relation for $\wed^kV$ in the form of
Hankel's determinant has an advantage of not using the elements that are not in the
ring of representations  of $\GL(V)$.

\begin{ex} \label{exgrothp1}
Let $\dim V=p|1$. Then ~\eqref{eqgrothrelation2} gives the relation
\begin{equation}\label{eqrelpone1}
    \begin{vmatrix}
      \wed^kV &   \wed^{k+1}V \\
      \wed^{k+1}V   & \wed^{k+2}V \\
    \end{vmatrix}=0\,,
\end{equation}
i.e., $\wed^kV \cdot\wed^{k+2}V=(\wed^{k+1}V)^2$ (product means tensor product) for
$k\geq p$. This can be seen directly as follows. $V=U\oplus \Pi W$ where $\dim U=p$,
$\dim W=1$. Hence $\wed^kV=\dsum\limits_{i+j=k} \wed^iU\otimes \Pi^j S^jW$. Note that
$S^jW=W^j$.  Thus for $k\geq p$ we have $\wed^kV=\dsum\limits_{i=0}^p \wed^iU\otimes
(\Pi W)^{k-i}$, therefore $\wed^{k+1}V=\wed^kV\otimes \Pi W$ (a geometric progression).
Obviously, by tensor multiplying $\wed^{k}V$ and $\wed^{k+2}V$  we get the isomorphisms
$\wed^{k}V\otimes\wed^{k+2}V=\wed^{k}V\otimes\wed^{k+1}V\otimes \Pi
W=\wed^{k+1}V\otimes\wed^{k+1}V$, which is exactly the relation~\eqref{eqrelpone1}.
\end{ex}

Let us obtain the expansion at infinity for the rational function $R_V(z)$. For this,
we shall rearrange the numerator and denominator in~\eqref{eqgrothfrac}. Since
$\wed^i(U)=\det U \otimes \wed^{p-i}(U^*)$ and $\wed^j(W)=\det W \otimes
\wed^{q-j}(W^*)$, we have
\begin{multline*}
    R_V(z)= \frac{\det U}{\det W}\,
    \frac{\wed^p(U^*)+z\wed^{p-1}(U^*)+ \ldots+z^p}{\wed^q(W^*)-z\Pi
    \wed^{q-1}(W^*)+ \ldots+(-z)^q\Pi^q}= \\
    \Ber V\,(-\Pi)^q\,z^{p-q}\,
    \frac{1+z^{-1}U^*+z^{-2}\wed^2U^*+\ldots+z^{-p}\wed^pU^*}{1-z^{-1}\Pi
    W^*+z^{-2}\wed^2W^*-\ldots+(-z)^{-q}\Pi^q\wed^qW^*}=\\
    \Ber
V\,(-\Pi)^q\,z^{p-q}\,\frac{\wed_{\frac{1}{z}}(U^*)}{\wed_{-\frac{1}{z}\Pi}(W^*)}=
    \Ber
    V\,(-\Pi)^q\,z^{p-q}\,\wed_{\frac{1}{z}}(V^*)=\\
    (-\Pi)^q\,z^{p-q}\,\Ber V \sum_{k\leq 0}z^k \wed^k (V^*)=
    (-\Pi)^q\, \sum_{k\leq 0} z^{k+p-q}\,\Sigma^{p-k} (V)=\\
    (-\Pi)^q\, \sum_{k\leq p-q} z^{k}\,\Sigma^{k+q} (V)\,.
\end{multline*}
Hence the rational function $R_V(z)$ taking values in a Grothendieck ring has  the
following expansions:
\begin{align}\label{eqexpansiongroth}
  R_V(z) & = \sum_{k\geq 0}\, z^k\,\wed^k(V) & \text{(at zero)}\\
    & = \sum_{k\leq p-q} z^{k}\,(-\Pi)^q\,\Sigma^{k+q}(V) &
    \text{(at infinity)}
\end{align}
In the same way as in Section~\ref{secber}  we arrive at the following theorem.

\begin{thm}\label{thmgroth2}
The sequence in the Grothendieck ring
\begin{equation}\label{eqdiffergroth}
    \G_k=\wed^k V-(-\Pi)^q\Sigma^{k+q} V
\end{equation}
for all $k\in \Z$ is a recurrent sequence of period $q$. \qed
\end{thm}

It very well fits with the equality $\wed^k V= \Sigma^kV$ of the classical case of
$q=0$, i.e., $\wed^k V= \det V\otimes \wed^{p-k} V^{*}$, which is a canonical
isomorphism compatible with the action of $\GL(V)$. Theorem~\ref{thmgroth2} implies
the vanishing of the Hankel determinants of order $q+1$ made of the elements $\G_{k}$.

\begin{ex} \label{exgroth11}
Consider  $V$ where $\dim V=1|1$. Then $\wed^k(V)=0$ for $k<0$, $\dim\wed^0(V)=1$,
$\dim\wed^k(V)=1+\Pi$ for $k\geq 1$. In the same way $\dim\Sigma^{k+1}(V)=1+\Pi$ for
$k\leq -1$, $\dim\Sigma^{1}(V)=1$, $\dim\Sigma^{k+1}(V)=0$ for $k>0$. It follows that
$\dim\wed^k V-(-\Pi)\dim\Sigma^{k+1}V=1+\Pi$ for all $k\in\Z$, which is a geometric
progression with ratio $\Pi$ infinite in both directions. This verifies the statement
of Theorem~\ref{thmgroth2} for $V$ at the    level of dimensions.
\end{ex}

\begin{ex} (Continuation of Examples~\ref{exgrothp1}
and~\ref{exgroth11}.) For a superspace $V$ such that $\dim V=1|1$ we shall show
explicitly an isomorphism $\f\co \wed^kV
\otimes\wed^{k+2}V\to\wed^{k+1}V\otimes \wed^{k+1}V$ commuting with the action
of $\GL(V)$. Let $e\in V_0$, $\e\in V_1$ be a basis of $V$. Then
$E_k=\underbrace{\e\wedge \ldots\wedge\e}_{k}$  and
$F_k=e\wedge\underbrace{\e\wedge \ldots\wedge\e}_{k-1}$  can be taken as a
basis in $\wed^k V$ for $k\geq 1$.   The desired isomorphism $\f$ can be
written as follows:
\begin{align*}
    \f(E_k\otimes E_{k+2}) &= \alpha\, E_{k+1}\otimes E_{k+1}\,,\\
    \f(E_k\otimes F_{k+2}) &=
\frac{1}{2}\,\Bigl(-\alpha+\frac{k}{k+1}\,\beta\Bigr)\,E_{k+1}\otimes
    F_{k+1} \\
    {\quad } & +(-1)^k
\frac{1}{2}\,\Bigl(\alpha+\frac{k}{k+1}\,\beta\Bigr)\,F_{k+1}\otimes
E_{k+1}\,,\\
    \f(F_k\otimes E_{k+2}) &=
(-1)^k\frac{1}{2}\,\Bigl(\alpha+\frac{k+2}{k+1}\,\beta\Bigr)\,E_{k+1}\otimes
    F_{k+1} \\
    {\quad } & +
\frac{1}{2}\,\Bigl(\alpha-\frac{k+2}{k+1}\,\beta\Bigr)\,F_{k+1}\otimes
E_{k+1}\,,\\
    \f(F_k\otimes F_{k+2}) &= \beta \,F_{k+1}\otimes F_{k+1} \,,
\end{align*}
where $\alpha,\beta$ are arbitrary nonzero parameters. In particular, notice
that $\f$ is not unique.
\end{ex}

\section{Cramer's Rule in Supermathematics} \label{seccramer}

In this section we formulate Cramer's rule in supermathematics basing on the
geometrical meaning of the Berezinian. Earlier such a generalization was
obtained by Bergveldt and Rabin in~\cite{rabin:duke}, who  used the `hard
tools' of the Gelfand--Retakh quasi-determinants theory
(see~\cite{gelfand:retakh91, gelfand:retakh02}).
 Our approach does not use anything but the main properties of the Berezinian.

Let us first formulate the usual Cramer's rule geometrically. Let $A$ be a linear
operator in an $n$-dimensional vector space $V$. Consider a linear equation
\begin{equation*}
A(\bx)=\by\,.
\end{equation*}
Here $\bx, \by$ are vectors in $V$. For any volume form $\r$ on $V$ and arbitrary
vectors $\bv_1, \ldots, \bv_{n-1}$ we obviously have
\begin{equation*}
\r(A(\bx),A(\bv_1),\ldots,A(\bv_{n-1}))=\det A\cdot \r(\bx,\bv_1,\ldots,\bv_{n-1}).
\end{equation*}
Considering this equation for different vectors $\bv_1, \ldots, \bv_{n-1}$ we can
express $\bx$ via $\by=A(\bx)$. Namely, let $\be_1,\ldots,\be_n$ be an arbitrary basis
in $V$. Take as $\r$ the coordinate volume form, i.e., $\r(\be_1,\ldots,\be_n)=1$ and
for any other vectors the value of $\r$  equals the determinant of the matrix
consisting of the corresponding coordinate row vectors. Then for the $k$-th coordinate
of $\bx$ we have
\begin{equation*}
x^k=\r(\be_1,\ldots,\bx,\ldots,\be_n)
\end{equation*}
($\bx$ stands at the $k$-th place), hence
\begin{equation*}
x^k=\frac{1}{\det A}\,\r(A(\be_1),\ldots,\by,\ldots,A(\be_n))= \frac{1}{\det
A}\,\begin{vmatrix}
  a_1{}^1 & \ldots & a_1{}^n  \\
  \dots &\dots &\dots  \\
  y^1 & \ldots & y^n  \\
   \dots &\dots &\dots\\
  a_n{}^1 & \ldots & a_n{}^n \\
\end{vmatrix},
\end{equation*}
where at the \RHS\  the coordinates of $\by$ replace the $k$-th row of the
matrix of the operator $A$. This is exactly Cramer's rule. Here we use row
vectors rather than columns because it is more convenient in the supercase.

These considerations can be  generalized to the supercase as follows.

Let $V$ be a $p|q$-dimensional linear superspace. Consider a volume form $\r$. Recall
that in the supercase a volume form is defined as a function on bases such that  a
change of basis  is equivalent to the multiplying  by the Berezinian of the transition
matrix. For example, a coordinate volume form associated with a basis
$\be_1,\ldots,\be_{p+q}$ where $\be_1,\ldots,\be_{p}$ are even vectors and
$\be_{p+1},\ldots,\be_{p+q}$ are odd vectors, on vectors $\bv_1,\ldots,\bv_{p+q}$ of
another basis of the same format equals the Berezinian
\begin{equation*}
\Ber\begin{pmatrix}
  v_1{}^1 & \dots & v_1{}^{p+q} \\
  \dots & \dots & \dots \\
  v_{p+q}{}^1 & \dots & v_{p+q}{}^{p+q} \\
\end{pmatrix}\,.
\end{equation*}
Here $\bv_i=v_i{}^j\be_j$. It follows that a volume form is linear in the first
$p$ arguments and hence can be extended by linearity to arbitrary vectors (the
last $q$ arguments must remain linearly independent odd vectors!). In
particular, it is possible to insert an odd vector into one of the first $p$
``even'' positions.

As above, for any volume form $\r$ on $V$ and  vectors $\bv_1, \ldots, \bv_{p+q-1}$ of
the appropriate parity we have
\begin{equation*}
\r(A(\bv_1),\ldots,A(\bx),\ldots,A(\bv_{p+q-1}))=\Ber A\cdot
\r(\bv_1,\ldots,\bx,\ldots,\bv_{p+q-1})\,,
\end{equation*}
where the vector $\bx$ stands at the one of the first $p$ ``even'' places. $A$ is
assumed to be an even invertible operator. This leads to a solution of a linear
equation
\begin{equation}\label{eqsuperlin}
A(\bx)=\by\,.
\end{equation}
in the superspace $V$ as follows. Take as $\r$ the coordinate volume form
associated with a basis $\be_1,\ldots,\be_{p+q}$. Then
$\r(\be_1,\ldots,\bx,\ldots,\be_{p+q})=x^k$, if $k=1,\ldots,p$. Hence the
formula for the first $p$ coordinates of $\bx$ corresponding to the even basis
vectors is exactly the same as in the classical case. For $k=1,\ldots,p$
\begin{equation}\label{eqcramer1}
    x^k=\frac{1}{\Ber A}\,\r(A(\be_1),\ldots,\by,\ldots,A(\be_{p+q}))=
\frac{1}{\Ber A}\,\D_k(A,\by)\,,
\end{equation}
where
\begin{equation}\label{eqcramer2}
\D_k(A,\by)=\Ber\begin{pmatrix}
  a_1{}^1 & \ldots & a_1{}^{p+q}  \\
  \dots &\dots &\dots  \\
  y^1 & \ldots & y^{p+q}  \\
   \dots &\dots &\dots\\
  a_n{}^1 & \ldots & a_n{}^n \\
\end{pmatrix},
\end{equation}
($\by$ inserted at the $k$-th ``even'' position). To obtain the last $q$ coordinates of
$\bx$ corresponding to the odd basis vectors $\be_{p+1},\ldots,\be_{p+q}$,  consider
the space $\Pi V$ with reversed parity. Let $\r^{\Pi}$ be the coordinate volume form on
$\Pi V$ corresponding to the basis
$$
\be_{p+1}\Pi,\ldots,\be_{p+q}\Pi,\be_1\Pi,\ldots,\be_{p}\Pi.
$$
Now we have
\begin{equation*}
\r^{\Pi}\bigl(\be_{p+1}\Pi,\ldots,\bx\Pi,\ldots,\be_{p+q}\Pi,\be_1\Pi,\ldots,\be_p\Pi\bigr)
=x^k
\end{equation*}
for $k=p+1,\ldots,p+q$. Introducing the notation
\begin{equation*}
\r^*(\bv_1,\ldots,\bv_{p+q}):=
\r^{\Pi}(\bv_{p+1}\Pi,\ldots,\bv_{p+q}\Pi,\bv_1\Pi,\ldots,\bv_{p}\Pi)
\end{equation*}
and
\begin{equation}
\Ber^*M:=\Ber M^{\Pi}
\end{equation}
for a matrix $M$, we can rewrite this as $x^k=\r^*(\be_1,\ldots,\bx,\ldots,\be_{p+q})$,
$k=p+1,\ldots, p+q$. Hence for $k=p+1,\ldots, p+q$
\begin{equation}\label{eqcramer3}
    x^k=\frac{1}{\Ber^* A}\,\r^*(A(\be_1),\ldots,\by,\ldots,A(\be_{p+q}))=
\frac{1}{\Ber^* A}\,\D_k^*(A,\by)\,,
\end{equation}
where
\begin{equation}\label{eqcramer4}
\D_k^*(A,\by)=\Ber^*\begin{pmatrix}
  a_1{}^1 & \ldots & a_1{}^{p+q}  \\
  \dots &\dots &\dots  \\
  y^1 & \ldots & y^{p+q}  \\
   \dots &\dots &\dots\\
  a_n{}^1 & \ldots & a_n{}^n \\
\end{pmatrix},
\end{equation}
($\by$ inserted at the $k$-th ``odd'' position).
Formulae~\eqref{eqcramer1}--\eqref{eqcramer4} give a complete solution of the
equation~\eqref{eqsuperlin}.  Recall that the matrix of a linear operator is
defined by the formula $A(\be_i)=a_i{}^j\be_j$. Hence
$A(\bx)=A(x^i\be_i)=x^ia_i{}^j\,\be_j$ if $A$ is even.

\begin{Rem}
For even invertible matrices the operation $\Ber^*$ is the same as $\Ber^{-1}$.
However, for matrices that are not invertible, $\Ber^*$ can make sense, taking
a nonzero nilpotent value,  while $\Ber$ and $\Ber^{-1}$ are not defined.
\end{Rem}

The ``super'' Cramer's formulae~\eqref{eqcramer1}--\eqref{eqcramer4} motivate
the following definition. Let $D_{ij}(A)$ denote the matrix obtained from an
even matrix $A$ by replacing all elements in the $i$-th row by zeros except for
the $j$-th element replaced by $1$. Notice that $D_{ij}(A)$ may be odd
depending on positions of the indices $i,j$.

\begin{de}
The \textit{$(i,j)$-th cofactor} or \textit{adjunct} of an even $p|q\times p|q$ matrix
$A$ is
\begin{equation}\label{eqadj}
    (\adj A)_{ij}:=\begin{cases}
    \Ber D_{ij}(A) &\text{ \quad when $i=1,\ldots, p$}\\
    \Ber^* D_{ij}(A) &\text{ \quad when $i=p+1,\ldots, p+q$}
    \end{cases}
\end{equation}
\end{de}

In the previous notation, $(\adj A)_{ij}=\D_i(A,\be_j)$ for $i=1\ldots,p$ and
$(\adj A)_{ij}=\D_i^*(A,\be_j)$ for $i=p+1\ldots,p+q$. Notice that this notion
is not symmetrical w.r.t. rows and columns, so it might better be called the
``right adjunct''. We have the following formulae for the entries of the
inverse matrix:

\begin{equation}
(A^{-1})_{ij}=\begin{cases}
    \dfrac{(\adj A)_{ji}}{\Ber A} &\text{ \quad when $j=1,\ldots,
    p$}\\[12pt]
    \dfrac{(\adj A)_{ji}}{\Ber^* A} &\text{ \quad when $j=p+1,\ldots, p+q$}
    \end{cases}
\end{equation}

\begin{ex}
Consider a $1|1\times 1|1$ even matrix
\begin{equation*}
A=\begin{pmatrix}
  a_1{}^1 & a_1{}^2 \\
  a_2{}^1 & a_2{}^2 \\
\end{pmatrix}=\begin{pmatrix}
  a & \beta \\
  \gamma & d \\
\end{pmatrix}\,.
\end{equation*}
Then by formulae~\eqref{eqadj} we get
\begin{align*}
    (\adj A)_{11}&=\Ber\begin{pmatrix}
      1 & 0 \\
      \gamma & d
    \end{pmatrix}  =\frac{1}{d} \\
      (\adj A)_{12}&=\Ber\begin{pmatrix}
      0 & 1 \\
      \gamma & d
    \end{pmatrix}  =-\frac{\gamma}{d} \\
    (\adj A)_{21}&=\Ber^*\begin{pmatrix}
      a & \beta \\
      1 & 0
    \end{pmatrix}=\Ber \begin{pmatrix}
      0 & 1 \\
      \beta & a
    \end{pmatrix}=-\frac{\beta}{a^2}
    \\
    (\adj A)_{22}&=\Ber^*\begin{pmatrix}
      a & \beta \\
      0 & 1
    \end{pmatrix}=\Ber \begin{pmatrix}
      1 & 0 \\
      \beta & a
    \end{pmatrix}=\frac{1}{a}
\end{align*}
Thus for the transpose adjunct matrix we have:
\begin{equation*}
B=\begin{pmatrix}
  \frac{1}{d} & -\frac{\beta}{a^2} \\
  -\frac{\gamma}{d^2} & \frac{1}{a} \\
\end{pmatrix},
\end{equation*}
and
\begin{equation*}
AB=\begin{pmatrix}
  a & \beta \\
  \gamma & d \\
\end{pmatrix}\begin{pmatrix}
  \frac{1}{d} & -\frac{\beta}{a^2} \\
  -\frac{\gamma}{d^2} & \frac{1}{a} \\
\end{pmatrix}=
\begin{pmatrix}
  \frac{a}{d}-\frac{\beta\gamma}{d^2} & 0 \\
  0 & \frac{d}{a}-\frac{\gamma\beta}{a^2} \\
\end{pmatrix}=
\begin{pmatrix}
  \Ber A & 0 \\
  0 & \Ber^*A \\
\end{pmatrix},
\end{equation*}
as expected.
\end{ex}

\begin{Rem}
A different approach to Cramer's rule was suggested
in~\cite{kantor-trishin:det}. They defined certain `relative determinants' of
$A$  polynomially depending on $A$ and considered  the `$\lambda$-solutions'
$\bx$ satisfying $A(\bx)=\lambda \cdot\by$ instead of $A(\bx)= \by$,  $\lambda$
being one of the relative determinants. This allowed them to avoid division and
to use only polynomial expressions.
\end{Rem}

\appendix
\section{Elementary Properties of Recurrent
Sequences} \label{append}

It is a classical result due to Kronecker that a power series represents a
rational function if and only if the infinite Hankel matrix  of the
coefficients has finite rank. In this section we summarize the relations
between  recurrent sequences and rational functions  used in the main text. We
present the material in the form convenient for our purposes.   Notice that
classical expositions (see~\cite{gantmacher:mat}) make use of the expansion of
a rational function at infinity, while we need to consider simultaneously two
expansions, at zero and at infinity.

Let
\begin{equation}\label{eqfraction}
    R(z)=\frac{a_0+a_1z+\ldots+a_pz^p}{b_0+b_1z+\ldots+b_qz^q}
\end{equation}
be a rational function. We assume that the numerator has degree $p$ and the denominator
degree $q$. The coefficients can be in an arbitrary commutative ring with unit.
Consider formal power expansions of the fraction~\eqref{eqfraction} at zero and at
infinity. Let $R(z)=\sum\limits_{k\geq 0}c_kz^k$ (near zero) and
$R(z)=\sum\limits_{k\leq p-q}c_k^*z^k$ (near infinity). Here and below it is convenient
to assume that coefficients such as $a_k$, $b_k$, $c_k$, etc., are defined for all
values of $k\in\Z$ but may be equal to zero for some $k$. Hence we have the equalities
\begin{equation}\label{eqanzero}
a_n=\sum_{i=0}^q b_i c_{n-i}
\end{equation}
for all $n$, where $c_k=0$ for $k<0$, and
\begin{equation}\label{eqaninf}
a_n=\sum_{i=0}^q b_i c_{n-i}^*
\end{equation}
for all $n$, where $c_k^*=0$ for $k>p-q$. Taking into account that $a_n=0$ for $n>p$ or
$n<0$, we obtain, respectively, that
\begin{equation*}
\sum_{i=0}^q b_i c_{n-i}=0
\end{equation*}
for all $n>p$, i.e.,
\begin{equation}\label{eqpqrelation}
\sum_{i=0}^q b_i c_{k+q-i}=0
\end{equation}
for all $k>p-q$, and that
\begin{equation}\label{eqpqrelation2}
 \sum_{i=0}^q b_i c_{k-i}^*=0
\end{equation}
for all $k<0$. Also, if we subtract~\eqref{eqaninf} from ~\eqref{eqanzero}, we obtain
that
\begin{equation}\label{eqreldif}
\sum_{i=0}^q b_i \dif_{n-i}=0
\end{equation}
for all $k\in\Z$, where $\dif_k=c_k-c_k^*$.

It is convenient to introduce the following definition. We say that a sequence
$\{c_k\}_{k\in\Z}$ is right or positive if $c_k=0$ for $k<0$.
\begin{de}
A right sequence $\{c_k\}$ is a \textit{$p|q$-recurrent sequence} or, shortly, a
\textit{$p|q$-sequence} if the elements $c_k$ satisfy a recurrence relation of the
form~\eqref{eqpqrelation} for all $k>p-q$.
\end{de}
It follows that the coefficients $c_k$ of the power expansion at zero of the
fraction~\eqref{eqfraction}  make a $p|q$-recurrent sequence. (The coefficients
of the expansion of~\eqref{eqfraction} at infinity also make a $p|q$-sequence
after the re-indexing that makes them a right sequence, $c_k':=c_{p-q-k}^*$.)
The fraction~\eqref{eqfraction} is classically  referred to as the
\textit{generating function} or the \textit{symbol} of the  recurrent sequence
$\{c_k\}$.

For a  sequence $\{c_k\}_{k\geq 0}$ to be a $p|q$-sequence means, if $p\geq q$, that it
satisfies a recurrence relation of period $q$ except for the $p-q+1$ initial  terms
$c_0, \ldots, c_{p-q}$, and if $p<q$, that it satisfies a recurrence relation of period
$q$ for all terms $c_k$, $k\geq 0$, and  can be extended to the left by $q-p-1$ zero
terms so that the relation   still holds. If we denote the set of all $p|q$-sequences
by $\seq_{p|q}$, then
\begin{equation*}
\seq_{p|q}\subset \seq_{p|q+1} \text{ and } \seq_{p|q}\subset \seq_{p+1|q} \,.
\end{equation*}

Hence we have the following picture for the coefficients of the expansions of
the rational function~\eqref{eqfraction}. The coefficients of the expansions at
zero and at infinity satisfy the same recurrence relations of period $q$. If
$p<q$, the coefficients $c_k$ and $c_k^*$ can be nonzero only in the disjoint
ranges $k\geq 0$ and $k\leq p-q$, respectively. The recurrence relation holds
for all terms. If $p\geq q$ (that is, when the fraction is improper), the
coefficients $c_k$ and $c_k^*$ can be simultaneously nonzero   in the finite
range $0\leq k\leq p-q$. Separate  recurrence relations   break down in this
range. However, in all cases the sequence $\dif_k=c_k-c_k^*$, infinite in both
directions and which  coincides with either $c_k$ or $-c_k^*$ `almost
everywhere', satisfies the recurrence relation for all $k\in\Z$.

If a sequence $\{c_k\}$ is given, one can consider the associated infinite
Hankel matrix with the entries $a_{ij}=c_{i+j}$. Let $\{c_k\}$ satisfy a
recurrence relation of the form~\eqref{eqpqrelation} for all $k\geq N$. Assume
that $b_0$ is invertible. Then  the infinite vector
$\c_{N+q}=\{c_{q+k}\}_{k\geq N}$ is a linear combination of the vectors
$\c_{N}=\{c_{k}\}_{k\geq N}$, \ldots, $\c_{N+q-1}=\{c_{q+k-1}\}_{k\geq N}$.
Hence their exterior product vanishes. In particular it implies the vanishing
of the Hankel minors of order $q+1$:
\begin{equation*}
\begin{vmatrix}
      c_k  & \dots & c_{k+q}  \\
      \dots & \dots & \dots \\
      c_{k+q}  & \dots & c_{k+2q}  \\
    \end{vmatrix}=0
\end{equation*}
where $k\geq N$.   On the other hand, solving a recurrence relation of period
$q$ involves division by a  Hankel determinant of order $q$. There is a vast
literature devoted to theoretical and practical aspects of recurrent sequences
and Hankel matrices.




\def\cprime{$'$} \def\cprime{$'$}

\end{document}